\documentclass[preprint]{imsart}

\RequirePackage{amsthm,amsmath,amsfonts,amssymb}
\RequirePackage[numbers]{natbib}

\RequirePackage{graphicx}

\usepackage{bbm}
\usepackage{aligned-overset}
\usepackage{color}
\usepackage{caption}
\usepackage{ulem}
\usepackage{subcaption}
\usepackage{booktabs}
\usepackage{todonotes}

\startlocaldefs
\theoremstyle{plain}
\newtheorem{theorem}{Theorem}

\newtheorem{lemma}{Lemma}
\newtheorem*{lemma_nn}{Lemma}
\newtheorem{corollary}{Corollary}

\theoremstyle{remark}
\newtheorem{step}{Step} 

\newtheorem{remark}{Remark}[]
\newcommand{\Acal}{\mathcal{A}}
\newcommand{\R}{\mathbb R}
\newcommand{\N}{\mathbb N}
\renewcommand{\L}{\mathcal{L}}
\newcommand{\Lmu}{\L_\#\mu}

\newcommand*\diff{\mathop{}\!\mathrm{d}}
\newcommand{\im}{\operatorname{im}}

\newcommand{\e}{\operatorname{e}}

\endlocaldefs

\begin{document}

\begin{frontmatter}
\title{Nested Sampling And Likelihood Plateaus}
\runtitle{Nested Sampling And Likelihood Plateaus}

\begin{aug}
\author[A]{\fnms{Doris} \snm{Schneider}\ead[label=e1]{doris.schneider@fau.de}},
\and
\author[A]{\fnms{Philipp} \snm{Wacker}\ead[label=e2, mark]{wacker@math.fau.de}}

\address[A]{Department of Mathematics, Friedrich-Alexander-Universität Erlangen-Nürnberg, Cauerstraße 11, 91058 Erlangen, Germany.
\printead{e1,e2}}

\end{aug}

\begin{abstract}
Nested sampling solves the problem of computing the evidence, also called the marginal likelihood, which is the integral of the likelihood with respect to the prior. The main idea of nested sampling is to replace the high-dimensional likelihood integral over the parameter space with an integral over the unit line by employing a push-forward with respect to a suitable transformation. For this transformation to work, it is often implicitly or explicitly assumed that samples from the prior are uniformly distributed on the unit line after having been mapped by this transformation. We show that this assumption is wrong if the likelihood function has a prior-non-negligible plateau. We show that the substitution enacted by nested sampling is valid nevertheless, thereby giving a rigorous justification of the nested sampling paradigm in a broader context than previously available. We also present a modification of nested sampling which improves its performance in the case of a likelihood function with a plateau.
\end{abstract}


\end{frontmatter}

\section{Introduction}
A recurring problem in computational statistics is evaluation of an integral of the form
\begin{equation}
    \mu(\phi) = \int_\Omega \phi(x) \diff\mu(x)
\end{equation}
where $(\Omega, \mathcal A, \mu)$ with $\Omega\subset \R^d$ is a probability space and  $\phi :\Omega\to \R$ is a measurable function. This is particularly challenging if regions of large measure $\mu$ and regions of interest for $\phi$ (i.e. peaks) do not coincide and generally in high dimensions $d \gg 1$.

We consider a more specific Bayesian setting where $\mu$ is interpreted as a prior probability measure and $\phi$ is a likelihood function $\L$. Note that $\L$ will in practice depend on some data but as we will assume $\L$ to be fixed, we suppress this dependence. Bayes' theorem states that the posterior measure $\mu_\text{post}$ is given by

\begin{equation}
    \mu_\text{post}(\diff x) = \frac{\L(x)\cdot \mu(\diff x)}{\int_\Omega \L(x)\mu(\diff x)}.
\end{equation}

The denominator $Z:=\int_\Omega \mathcal{L}(x) \diff \mu(x)$ is usually called the marginal likelihood or  ``evidence'' because of its interpretability as the \textit{evidence of the data for a given model} in the context of Bayesian model selection. While many sampling methods do not require $Z$ to be known, it is needed for some tasks like Bayesian model selection, see for example \cite{mackay2003information} for a lucid exposition.

Nested sampling was introduced as a specialized algorithm for the computation of the evidence in \cite{skilling2006nested} and has been applied succesfully in astronomy and computational physics (\cite{mukherjee2006anested, feroz2009multinest, feroz2008multimodal, parkinson2013bayesian,vegetti2009bayesian,veitch2015parameter,baldock2017constant,murray2006nested,veitch2015parameter,partay2014nested}), biomathematics (\cite{aitken2013nested, dybowski2013nested, pullen2014bayesian}) and other fields. It seems to not yet have gained much traction in the statistical community, with \cite{salomone2018unbiased} giving some reasons for this fact. From a bird's eye view, nested sampling is a method of computing the evidence $\int_\Omega \L(x)\diff \mu(x)$ by converting this high-dimensional integral into an integral on the unit line $[0,1]$:
\begin{equation} \label{eq:Nested-sampling-steps}
     \int_\Omega \mathcal{L}(x) \diff \mu(x) = \int_0^1 \tilde{\mathcal{L}}(X) \diff X \approx  \text{Monte Carlo approximation.} \end{equation}
Here, $ \tilde \L$ is a reparametrized version of $\L$ defined below. In \cite{skilling2006nested}, $\tilde\L$ is the ``overloaded'' form of the likelihood $\L$ but here we will explicitly notationally distinguish between $\L$ and $\tilde \L$.
 
There is already an extensive body of work regarding the performance of the second step, i.e.\ the quality of the Monte Carlo approximation of the 1d integral $\int_0^1 \tilde \L(X) \diff X$ and a series of improvements to Skilling's original method has emerged, for example in \cite{higson2018sampling,feroz2013importance,higson2019dynamic,salomone2018unbiased}. To the best of our knowledge there has not been much discussion on the validity of the first step, i.e.\ the substitution of the integral over $\Omega$ with an integral over the unit interval. \cite{salomone2018unbiased} gives a very short and precise derivation via the inverse of the survival function but does not consider the problematic case of discontinuous survival functions (corresponding to plateaus in the likelihood, see below). In \cite{chopin2010properties}, the authors refer to \cite{burrows1980new} for justification of the integral transformation, but in our opinion the issue is quite difficult and deserves a more careful analysis. At the heart of the matter, nested sampling amounts to integration with respect to the push-forward measure, but the details are surprisingly involved and contain many steps that ``almost'' fail to hold. 

The main difficulty arises when $\L$ has a non-negligible plateau, i.e. if there exists a level $\alpha\in \R$ such that $\mu(\L = \alpha) > 0$. This was already hinted at (although in the context of the algorithm's performance) in the original publication \cite{skilling2006nested} where they discuss ``cliffs'' -- deemed non-problematic -- and ``plateaus''. The authors recognize the difficulty of plateaus (for reasons laid out below), but they rule ``[...] even so, it may be possible to generate [new active samples] efficiently.'' In other words, they state that the existence of plateaus makes the Monte Carlo method more prone to performance problems which can be overcome by a cleverer method. Murray's PhD thesis \cite{murray2007advances} also mentions this issue and further improves on Skilling's original suggestion.

We will show that plateaus in $\tilde{\L}$ are not only computationally troublesome but also fundamentally mathematically problematic: If there is a level $\alpha \in \R$ such that $\mu(\L = \alpha) > 0$, then the following implicit uniformity assumption about nested sampling is violated:

``Samples from the prior which are plotted in a $X$-$\L$-diagram are uniformly distributed along the axis $[0,1]$.'' (The meaning of $X$ is explained in detail later). This is what Sivia and Skilling mean when they write ``In terms of $\xi$, the objects are uniformly sampled subject to the constraint $\xi < \xi^\star$'' \cite[section 9.2]{sivia2006data}. 

We show that this is indeed wrong in general. As the assumption is critical to the original justification of the correctness of nested sampling, this seems to invalidate the main paradigm of nested sampling (at least in the case of problematic plateauing likelihood functions). Fortunately, we can prove that the integral substitution which is the first equality in \eqref{eq:Nested-sampling-steps} is still correct. 

But, although this resolves theoretical objections to the applying nested sampling in the case of likelihood functions with a plateau, there are practical problems. This has been described \cite{skilling2006nested} and \cite{murray2007advances}, who suggest to resolve with a randomization or ``labelling'' approach, breaking a tie between points with the same likelihood value. We argue in this manuscript that there is a computationally more suitable way by splitting the integration domain into plateaus and their complement. Both components can then be evaluated separately in more efficient manner: The plateau component consists of the integration of constant functions which is done in a trivial manner, the non-plateau component can be computed by either a vanilla nested sampling method or by one of the various improvements like MultiNEST \cite{feroz2009multinest}, NS-SMC \cite{salomone2018unbiased}, dynamic nested sampling \cite{higson2019dynamic}, or others. This means that we are not proposing an alternative to nested sampling, but rather a preprocessing step in order to handle the likelihood plateaus. After that, the application nested sampling can (and should) be augmented with one of the algorithmic improvements mentioned. For example, \cite{salomone2018unbiased}, which derives a very appealing connection between nested sampling and the sequential Monte Carlo approach, propose a modification of nested sampling, NS-SMC, which is competitive with gold-standard Monte Carlo methods. Nevertheless, in their exposition (section 5.2) they assume exclusion of ``degenerate likelihoods'' by which they mean plateaus in the likelihood, i.e. $\mu(\L(x) = l) = 0$ for all levels $l\in \R$. Our proposed splitting preprocessing step creates exactly this setting and thus NS-SMC (or other competitive nested sampling variants) can be applied directly.

The contributions of this manuscript are, stated succinctly, as follows:

\begin{enumerate}
    \item We show that the integral transformation $\int_\Omega \mathcal{L}(x) \diff \mu(x) = \int_0^1 \tilde{\mathcal{L}}(X) \diff X$ at the heart of nested sampling is valid even for very general likelihood functions with plateaus, although this is not entirely straightforward to see.
    \item If the likelihood in question has a plateau of positive measure, we propose to split into plateau and non-plateau components first. We demonstrate that this improves the computational efficiency of nested sampling by a reasonable amount, especially in settings where evaluation of the likelihood is computationally expensive (for example if it entails numerical solution of a PDE).
\end{enumerate}

The remainder of the manuscript is structured as follows. After a quick motivation for why we want to consider likelihoods with plateaus at all, in section \ref{sec:nestedsampling} we give a quick recap of nested sampling, detailing how the main paradigm of nested sampling consists of a specific integral transformation whose correctness -- we argue -- needs to be proven. After that, section \ref{sec:mainParadigm} lays out the two main issues with this integral transformation, as handled by lemmas \ref{lem:integral_wrt_pshfwd} and \ref{lem:characterization-problematic} and states the main theoretical result, theorem \ref{thm:main-theorem}. Section \ref{sec:splitting} describes a pre-processing step suitable for settings with likelihood plateaus and demonstrates that our modification improves on the randomization procedure as proposed by \cite{skilling2006nested} and \cite{murray2007advances}. The appendix deals with the proofs for lemmas \ref{lem:integral_wrt_pshfwd} and  \ref{lem:characterization-problematic}  and thus for theorem \ref{thm:main-theorem}.

We suspect that the techniques employed (especially Lemma~\ref{lem:main}) have some connections with the notion of generalized inverse distribution functions (\cite{de2015study,kampke2015generalized,embrechts2013note}) which is slightly (but unintentionally) obfuscated by the fact that we consider non-increasing instead of non-decreasing functions. We haven't pursued this direction but it might be interesting to draw from this connection.

After completion of this manuscript, there has been more work \cite{fowlie2020nested} regarding the implementation of nested sampling in the case of likelihoods with plateaus with a similar way suggestion how to handle this case in practice.

\subsection{Likelihoods with plateaus}
While likelihood functions with plateaus may seem like an exotic special case, there are applications where they are the most fitting mathematical model. Essentially, there are two important scenarios where the likelihood function has a plateau of positive measure:

In the context of quantization in electrical engineering, signals are rounded to a grid. The likelihood of the true signal given such a quantized version is then a uniform distribution over the range of numbers being rounded to this value \cite{sripad1977necessary,bjorsell2007truncated}. This amounts to a likelihood function with a flat plateau as its peak (i.e. a set of positive measure is ``equally likely'' in the light of the data).

Second, some applications come with explicit knowledge about upper bounds on the measurement error's magnitude which corresponds to a plateau (of magnitude 0) on the range of impossible original signals (e.g. in biological applications; in the context of image processing on grayscale images with values in the unit interval; or when reading out an analog thermometer by looking at the nearest labeled tick on the scale: Here, the maximum measurement error is the distance between adjacent ticks on the scale). In \cite{Yao2017} the authors argue that bounded noise is more realistic in some biological and physical context and that the choice of the correct noise model has a large influence on the long-term behaviour of models. See also \cite{d2013bounded} and the references therein for more examples in physics, biology, and engineering. In this case the plateau arises by the likelihood's compact support.

\begin{figure}
    \centering
    \includegraphics[width=0.8\textwidth]{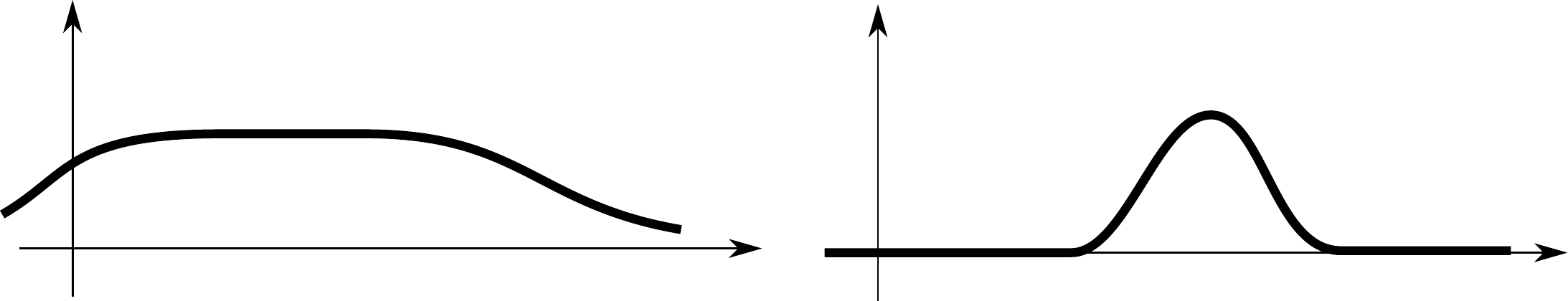}
    \caption{Left: Plateau as a peak. Right: Plateau by compact support.}
    \label{fig:plateaus}
\end{figure}

We also argue that, by treating the problematic case of likelihood functions with plateaus, we can close the previous gap in analytical verification of nested sampling's correctness.


\section{Nested sampling}\label{sec:nestedsampling}
Let $(\Omega, \mathcal A, \mu)$ be a probability space where we interpret $\mu$ as a prior probability measure of parameters $x\in\Omega\subset \R^d$. Sometimes (but not for our main results) we will assume that $\mu$ has a Lebesgue density $\pi$. Furthermore, we consider a likelihood $\L:\Omega \to \R$. This is usually due to the arrival of measurement data giving implicit information about possible parameter values. Hence, by Bayes' Law, the posterior probability is then given by $\diff \mu_\text{post}(x) = Z^{-1}\cdot \L(x) \cdot \diff \mu(x)$. The normalization constant $Z = \int_\Omega \L(x) \diff \mu(x)$ is called the evidence and is used for example in Bayesian model selection. It is usually a non-trivial matter to compute $Z$. Nested sampling is a method specifically designed for doing this.

\begin{figure}
        \centering
       \includegraphics[width=0.5\textwidth]{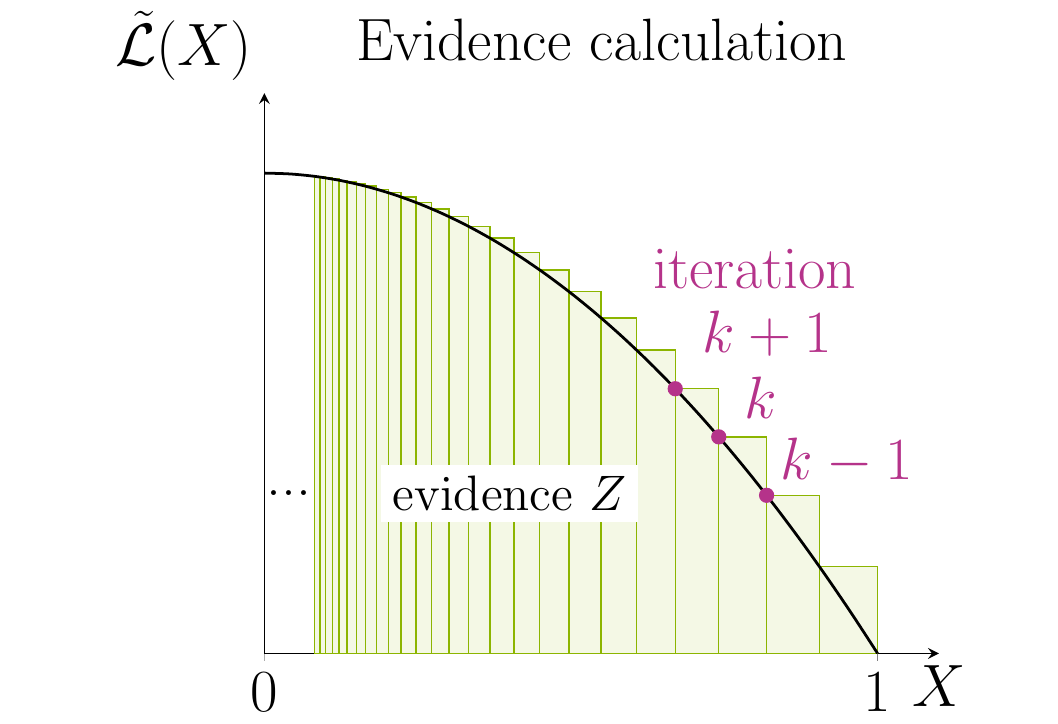}
        \caption{The integral from equation~\eqref{eq:paradigm} is approximated by an upper sum.}
        \label{fig:visualization_EvidenceUpperSum}
    \end{figure}

The basic idea of nested sampling is to swap the integration domain from $\Omega$ for the interval $[0,1]$ by constructing a specific transformation $\Phi$ and conducting basically a change of variables (or rather, a push-forward). Then the high-dimensional integration over the likelihood (which constitutes the evidence) becomes an integral of a function over the line $[0,1]$, which can more readily be approximated by a Monte-Carlo method. Concretely, nested sampling evaluates $Z$ by approximating the right-hand-side of 
\begin{equation} \label{eq:paradigm} 
    \int_\Omega \mathcal{L}(x) \diff \mu(x) = \int_0^1 \tilde{\mathcal{L}}(X) \diff X\end{equation}
instead of the left-hand-side (which is the definition of $Z$). Then this one-dimensional integral is efficiently approximated by a Monte-Carlo method (using ``active'' and ``dead'' samples and a clever way of estimating probabilities). An in-depth explanation of this remarkable idea can be found in for example \cite{skilling2006nested,sivia2006data}  but for sake of readability we record a quick outline of the nested sampling algorithm:
In the first step, $n \in \mathbb{N}$ samples $s_i, i = 1,...,n$ from the prior, the initial set of the so-called active samples, are generated and their corresponding likelihood values $\L(s_i)$ are calculated. In iteration $k$, the active sample with the lowest likelihood value $\L(\tilde{s}_k)$, is identified. Its accessible prior volume, which is a proportion of the prior given by $X_k = \int_{\L(x) > \L(\tilde{s}_k)} \diff\mu(x)$, is approximated by $X_k \approx \e^{-k/n}$. This approximation is directly based on the strong assumption (which we will show to be false in general) that $\Phi_\#\mu$ (introduced below) is the uniform measure on $[0,1]$. The sample $\tilde{s}_k$ is then replaced by a ``better'' sample $s_\mathrm{new}$ from the prior conditioned on $\L(s_\mathrm{new}) > \L(\tilde{s}_k)$. The number of active samples is again $n$ and the algorithm resumes with iteration $k+1$. The removed samples are called dead samples and are used for the computation of the evidence via $Z\approx \sum_k  (X_{k-1} - X_k) \L(\tilde{s}_k) = \sum_k (X_{k-1} - X_k) \tilde\L(X_k)$ as an upper sum of the integral in equation~\eqref{eq:paradigm} visualized in Fig.~ \ref{fig:visualization_EvidenceUpperSum}.


This figure jointly with figure \ref{fig:visualization-Xx} shows also the main idea of nested sampling: Folding (or nesting) all points on a level set of constant $\L$ into one point on the interval $[0,1]$ and hereby reducing dimensional complexity. We will concern ourselves with this transformation of integration and we show specifically under which constraints this construction works (regardless of the Monte-Carlo approximation procedure used for evaluation of the $1d$ integral).

We want to clarify once more that this manuscript is not primarily about convergence of the Monte-Carlo method involved with evaluating the integral $\int_0^1 \tilde\L(X)\diff X$, but we are interested in proving the equality in \eqref{eq:paradigm}.

\section{Correctness of nested sampling in the case of likelihood plateaus}\label{sec:mainParadigm}

We demonstrate next how we need to choose the transformation and the function $\tilde \L$ such that identity~\eqref{eq:paradigm} has a chance to hold.

\begin{figure}
        \centering
        \includegraphics[width=1\textwidth]{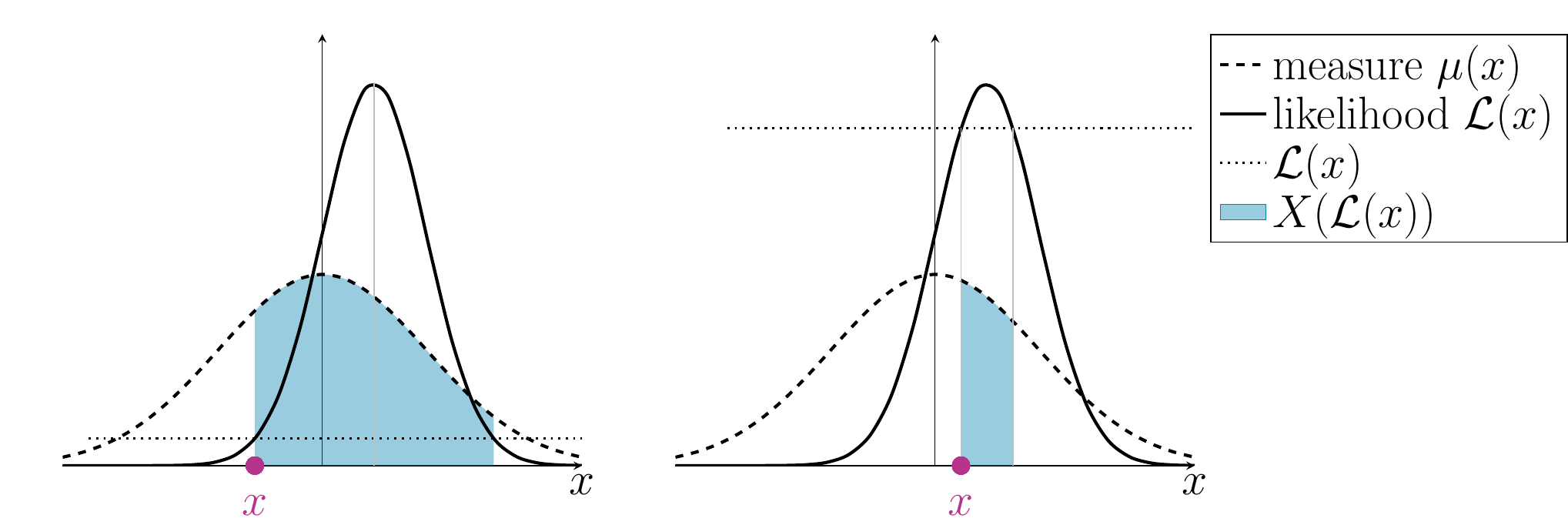}
        \caption{Example of the meaning of $\Phi(x) = X(\L(x))$ for a Gaussian measure $\mu$ and a Gaussian likelihood function $\L$ with $\Omega = \mathbb{R}$ for two representative $x\in \Omega$. }
        \label{fig:visualization-Xx}
    \end{figure}

First, we define $X(\lambda) = \mu(\{z\in\Omega: \L(z) > \lambda\})$, i.e.\ the $\mu$-measure of the $\lambda$-super-level-sets (also called the survival function in some contexts). Then we introduce a mapping $\Phi : (\Omega, \mathcal{B}^d) \to ([0,1], \mathcal{B}([0,1]))$ with	
	\begin{align} \notag
		\Phi(x) &= X(\mathcal{L}(x))\\ \notag
		&= \mu\left( \left\lbrace z \in \Omega : \mathcal{L}(z) > \mathcal{L}(x) \right\rbrace  \right)
	\end{align}
    which is visualized in Fig.~\ref{fig:visualization-Xx} for $\Omega = \mathbb{R}$, a Gaussian measure $\mu$ and a Gaussian-type likelihood. We remark that $\Omega = \mathbb{R}$ is chosen just for visualization reasons and all considerations of the proof are also valid in higher dimensions.
    
	\begin{figure}
	\centering
		\includegraphics[width=.9\textwidth]{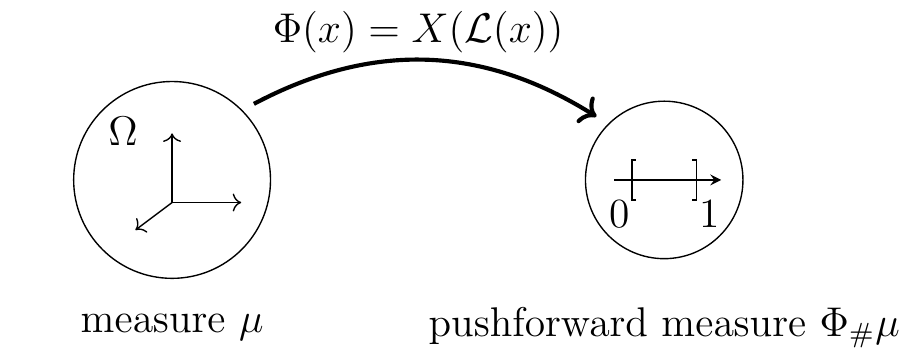}
		\caption{Visualization of relation of probability measure $\mu$ on $\Omega$ and the pushforward measure $\Phi_\#\mu$ on $[0,1]$. We call the idea of transforming the evidence integral on $\Omega$ into an integral on $[0,1]$ the main paradigm of nested sampling.} \label{fig:measure_mu_phimu}
	\end{figure}    
	
    This mapping $\Phi$ is the transformation which allows us to shift the integration from $\Omega$ to $[0,1]$ (see Fig.~\ref{fig:measure_mu_phimu} for a visualization), by virtue of the push-forward measure $\Phi_\#\mu := \mu \circ \Phi^{-1}$, a probability measure on $([0,1], \mathcal B([0,1]))$. 
    
        Next we need to specify the form of $\tilde \L$, the new function which is integrated over the line $[0,1]$. Skilling also calls it $\L$ (an ``overloaded form'') and identifies it with the inverse function of $X\circ \L$, but we need to be careful: In general, $X\circ \L$ has no inverse. The correct object to consider here is
    \begin{equation*}
    \tilde \L(\xi) = \sup\{\lambda \in \im \L: X(\lambda) > \xi\}.
    \end{equation*}
    $\tilde \L$ is a generalized inverse (\cite{chopin2010properties,kampke2015generalized,de2015study}) for $X$. We will show the following lemma:
    \begin{lemma_nn}[formal version of lemma \ref{lem:integral_wrt_pshfwd}]\label{lem:integral_wrt_pshfwd_formal}
    \begin{equation}\label{eq:integralequality}
        \int_0^1 \tilde \L(t) \diff \Phi_\#\mu(t) = \int_0^1 \tilde \L(t) \diff t.
    \end{equation}
\end{lemma_nn}
    It is not true that $\Phi_\#\mu$ is the uniform measure on $[0,1]$ (this is only the case for likelihood functions without plateau, see \cite[section 9.2]{sivia2006data}), which would in particular imply equation \eqref{eq:integralequality} more directly. We also prove in  
    \begin{lemma_nn}[formal version of lemma \ref{lem:characterization-problematic} in section \ref{sec:inverseX}]
    \label{lem:characterization-problematic_formal}
    For $\mu$-almost-all $x\in\Omega$, 
    \[\tilde \L(X(\L(x))) = \L(x). \]
    \end{lemma_nn}
    
    In the general case we cannot hope for this identity to hold everywhere (i.e.\ for all $x \in \Omega$) and then of course $\tilde \L(X(\lambda)) \neq \lambda$ for general $\lambda \in \R$. 
     
     Now we are able to derive the integral transformation which lies at the heart of nested sampling.

	\begin{align} 
	\int_\Omega \mathcal{L}(x) \diff \mu(x)  &{=}\int_\Omega\tilde{\mathcal{L}}\left( X(\mathcal{L}(x))\right) \diff \mu(x) \label{eq:line1}\\
	    &= \int_\Omega (\tilde{\mathcal{L}} \circ \Phi)(x) \diff \mu(x)\notag\\
	    &=\int_{[0,1]} \tilde{\mathcal{L}}(r) \diff (\Phi_\#\mu)(r)\notag\\
	    &= \int_0^1\tilde{\mathcal{L}}(X) \diff X
	      \label{eq:proofMainParadigm}
	\end{align}
	Equation~\eqref{eq:line1} holds because $\tilde\L(X(\L(x))) = \L(x)$ is true for $\mu$-almost-all $x\in\Omega$ (but not for all $x\in \Omega)$. This is proven in lemma \ref{lem:characterization-problematic}.
	The equality in \eqref{eq:proofMainParadigm} is shown to be true in lemma \ref{lem:integral_wrt_pshfwd} (Section~\ref{sec:uniform}) under quite mild assumptions, but not for obvious reasons: We will show that in general, $\Phi_\#\mu \neq \operatorname{Unif}[0,1]$ (which would in particular imply the equality).

\begin{figure}
        \centering
        \includegraphics[width=1\textwidth]{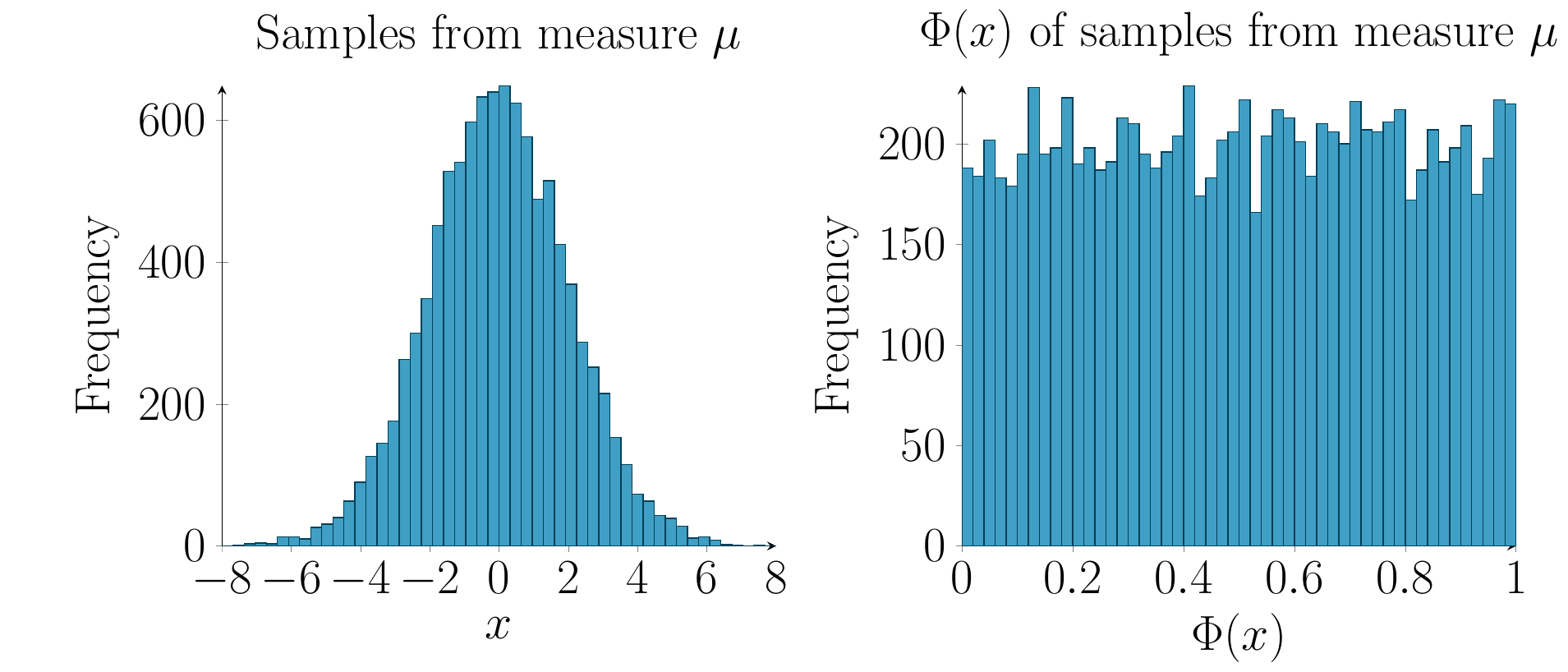}
        \caption{Histogram of 10000 samples from a Gaussian measure $\mu$ (left) and their map via $\Phi$ with Gaussian likelihood (mean 0 and variance 1) (right). The shape of the histogram corresponds to the shape of the probability densities of the Gaussian measure $\mu$ and the uniform pushforward measure $\Phi_\#\mu$, respectively. In the case of a likelihood with a plateau, the uniformity assumption on $\Phi_\#\mu$ is violated.}
        \label{fig:SamplesPriorPhi}
    \end{figure}

	The last step \eqref{eq:proofMainParadigm} can be interpreted in another way: 
	$\Phi$ is usually assumed (\cite{skilling2006nested,sivia2006data}) to have the following desirable property: If we generate i.i.d. samples $x_i\in \Omega$ from $\mu$ and map them via $\Phi$ into $[0,1]$, then the set $\{\Phi(x_i)\}_i$ consists of samples from a \textit{uniform} distribution on the set $[0,1]$ (see also Fig.~\ref{fig:SamplesPriorPhi}). This is equivalent to saying that the measure $\Phi_\#\mu$ is a uniform measure on $[0,1]$, or that the cumulative distribution $\Phi_\#\mu([0,r]) = r$. This, in turn, means that the Lebesgue-density of $\Phi_\#\mu$ is $\diff \Phi_\#\mu = \diff r$. In a consequence, we could also prove step \eqref{eq:proofMainParadigm} if $\Phi$ is a uniform measure. This is what Sivia and Skilling assume to hold when they write ``In terms of $\xi$, the objects are uniformly sampled subject to the constraint $\xi < \xi^\star$'' \cite[section 9.2]{sivia2006data}. 
	
	One of the main point of this manuscript is that this is sometimes wrong but that the integral transformation is still valid. More concretely, $ \Phi_\#\mu  \neq \text{Unif}[0,1]$, but it is a linear combination of uniform measures and Dirac deltas (i.e.\ point mass distributions).
	
	In one paragraph, the theoretical results of this manuscript can be summarized as follows:
	
	If $\L$ has a plateau of positive measure $\mu$, then $\Phi_\#\mu$ is not the uniform measure on $[0,1]$. Also, $\tilde\L(X(\L(x))) \neq \L(x)$ in general. Still, \eqref{eq:line1} and \eqref{eq:proofMainParadigm} are true and thus the main paradigm of nested sampling \eqref{eq:paradigm} holds.
	
    This constitutes our main theoretical result:

\begin{theorem}\label{thm:main-theorem}
    Let $(\Omega, \mathcal{B}^d, \mu)$ be a probability space with $\Omega \subseteq \mathbb{R}^d$, $d \in \mathbb{N}$ and $\mathcal{B}^d$ the Borel $\sigma$-algebra on $\Omega$. Let further $\mathcal{L} : (\Omega, \mathcal{B}^d) \to (\mathbb{R}, \mathcal{B})$ be a measurable function bounded from above and below, whereby $\mathcal{B}$ is the Borel $\sigma$-algebra on $\mathbb{R}$. Then
\begin{align}
\int_\Omega \mathcal{L}(x) \diff \mu(x) = \int_0^1 \tilde{\mathcal{L}}(X) \diff X.
\end{align}

\end{theorem}

\begin{proof}
	Follows from \eqref{eq:line1} with Lemma~\ref{lem:characterization-problematic} (in Section~\ref{sec:inverse}) and from \eqref{eq:proofMainParadigm} with Corollary~\ref{cor:inversion} (in Section~\ref{sec:uniform}).
\end{proof}

\section{Handling likelihood plateaus practically}\label{sec:splitting}
In this section we describe an easily implemented preprocessing step which transforms an integration problem of type $\int_\Omega \L(x)\diff \mu(x)$ with a plateau in $\L$ into one without a plateau (and a trivial additional term taking care of the plateau). 

We start by demonstrating that ``ignoring the plateau'' (i.e. running nested sampling without special consideration for the plateau) does not work. 

After that we show that the randomization strategy proposed by \cite{skilling2006nested} does indeed lead to better results but that it is computationally more expensive, as measured by the amount of likelihood evaluations needed. In models where the likelihood evaluation is costly (for example if model computation entails numerical solution of a partial differential equation), this leads to a high computational overhead, the reason being that a large fraction of sample candidates are discarded after evaluating their likelihood score. We then show how an additional preprocessing step leads to strongly improved performance.

\subsection{Splitting the parameter space}
\begin{figure}
    \centering
    \includegraphics[width=0.5\textwidth]{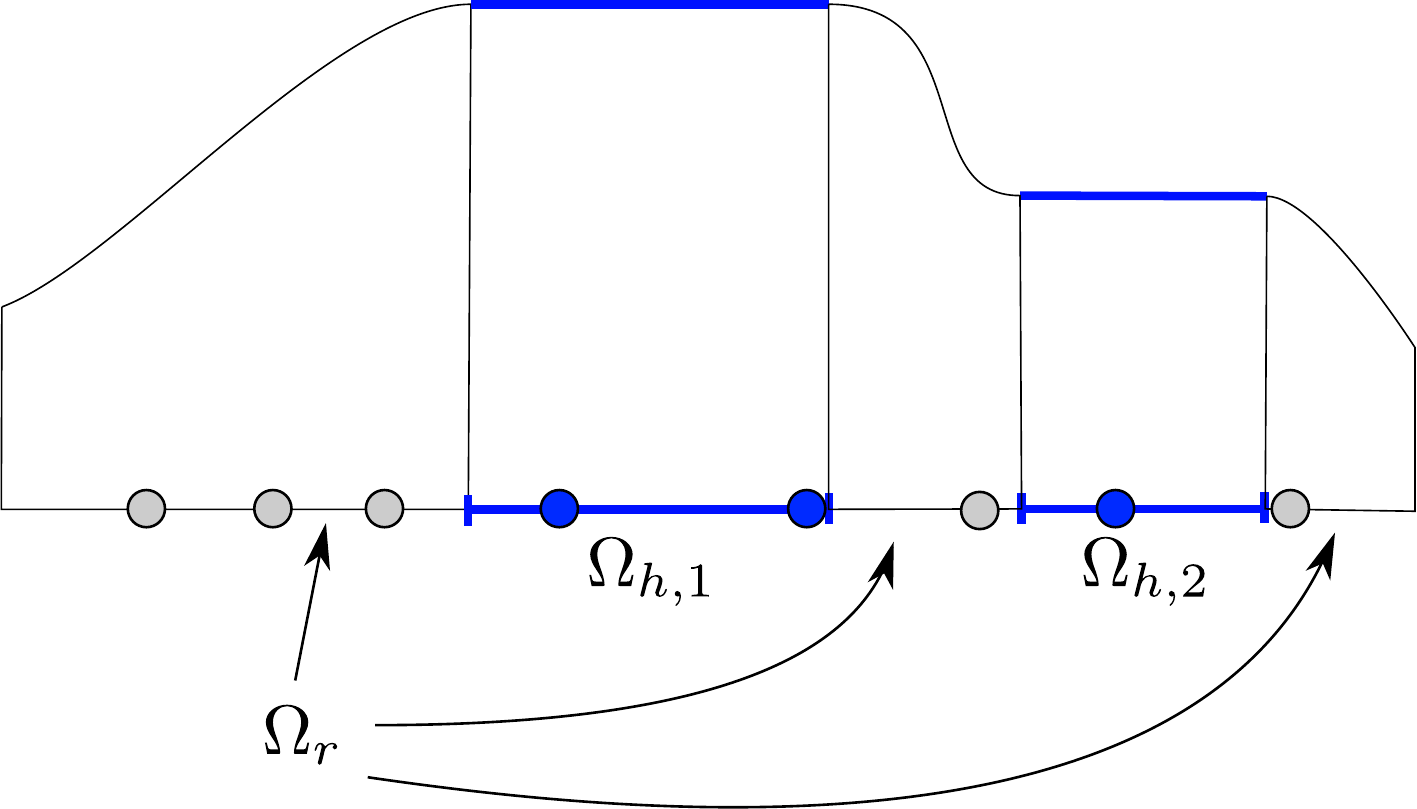}
    \caption{The $M$ initially generated samples (circles) are sorted into two categories: Those lying in one of the plateau domains $\Omega_{h,c}$ (blue circles) are used for the estimation of $\Delta_h$. The rest (grey circles) becomes the active set for an implementation of nested sampling on the non-plateau domain $\Omega_r$. The graph shown is the graph of the likelihood $\L$. The prior is not visualized here.}
    \label{fig:splitting}
\end{figure}

The main idea of the preprocessing step is to decompose the parameter set $\Omega$ into two disjoint subsets: on $\Omega_c \subseteq \Omega$ (c for ``constant'') the likelihood function has plateaus and on $\Omega_r = \Omega \setminus \Omega_c$ (r for ``regular'') it does not. The decomposition of $\Omega$ is illustrated in Figure~\ref{fig:splitting}. The computation of the evidence is conducted separately on either subset and merged subsequently. 

\begin{figure}
    \centering
    \includegraphics[width=0.5\textwidth]{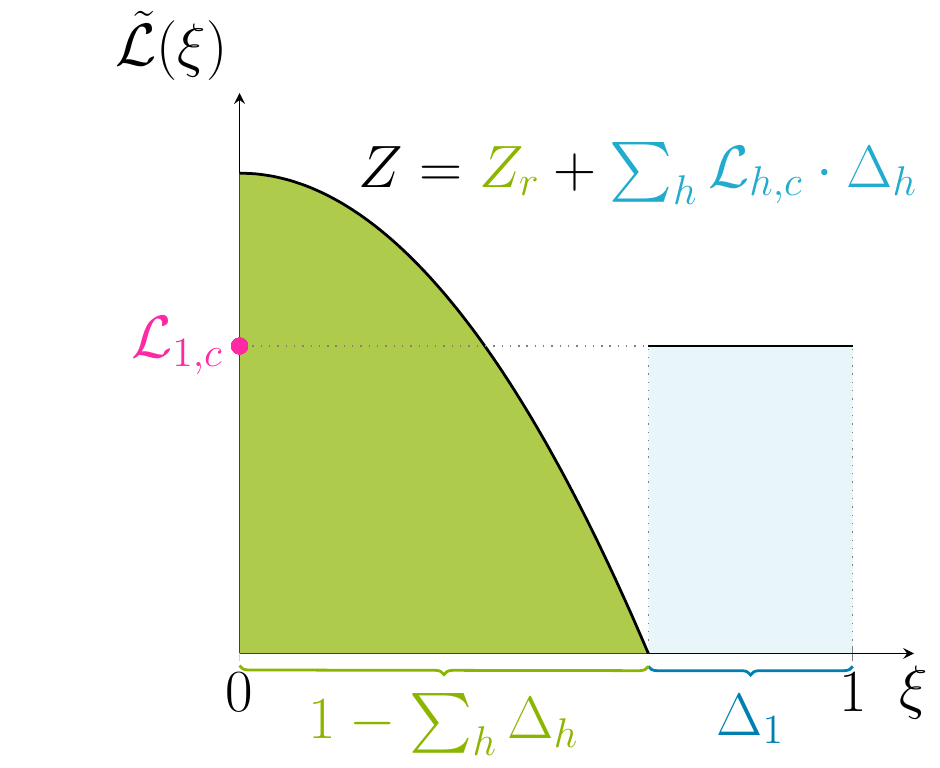}
    \caption{Illustration of evidence calculation with the splitting preprocessing for a likelihood with a plateau.}
    \label{fig:decompositionEvidence}
\end{figure}

On $\Omega_c$, integration should be easy because the most difficult integrand (the likelihood) is piece-wise constant, although in practice we do not know the shape and measure of $\Omega_c$ which usually precludes analytical computation. Hence, we further decompose $\Omega_c = \bigcup_{h=1}^H \Omega_{h,c}$ with $\Omega_{h,c} := \{x\in \Omega: \L(x) = \L_{h,c}\}$ assuming a finite number of likelihood plateau heights $h = 1, ..., H, H \in \mathbb{N}$ denoted by $\L_{h,c}$. Note that $r_i = \L_{i,c}$ in Lemma~\ref{lem:main} but we opted for the slightly more expressive notation $\L_{i,c}$ for $i=1,\ldots, H$ here. We additionally denote the prior volume of the $h$-th plateau component by $\Delta_{h} = \mu(\Omega_{h,c})$ resulting in the following integral decomposition also visualized in Figure~\ref{fig:decompositionEvidence}
\begin{align}\notag
    Z &= \int_\Omega \L(x) \diff \mu(x)\\ \notag
        &= \int_{\Omega_r} \L(x) \diff \mu(x) + \int_{\Omega_c} \L(x) \diff \mu(x)\\
        &\approx Z_r + \sum_h \L_{h,c} \cdot \Delta_{h} \label{eq:robustImplementation_calculationEvidence}
\end{align}

On $\Omega_r$, any version of nested sampling (either vanilla or a computationally improved modification like MultiNEST \cite{feroz2009multinest}, NS-SMC \cite{salomone2018unbiased}, or dynamic nested sampling \cite{higson2019dynamic}) can be performed to obtain $Z_r$ with the only necessary modification being the restriction of the initial accessible prior mass of $\Omega_r$ given by $1 - \sum_h \Delta_h$ (instead of 1). Nevertheless, the accessible prior volume keeps shrinking by the factor $\e^{1/n}$ from iteration to iteration since the uniformity assumption holds by Lemma~\ref{lem:main}. 

Finally, the evidence calculation in the presence of likelihood plateaus reduces to the efficient estimation of the plateau weights $\Delta_{h,c}$. Here, we propose to use relative frequencies as estimators because they are side products from the generation of the initial active set for nested sampling on $\Omega_r$. We draw $m$ samples according to $\mu$ and split the samples into groups: there are $n_h$ samples in $\Omega_{h,c}$ leading to the relative frequencies $\Delta_{h,c} = \frac{n_h}{m}$. The remaining $n_r$ samples are all in $\Omega_r$ and, therefore, are an initial active set for performing nested sampling. This automatically weights the accessible prior mass of $\Omega_r$: the higher $n_r$, the higher the accessible prior mass, the less uncertain the nested sampling algorithm.

Note again that we are not proposing an alternative to recent and competitive modifications to nested sampling (MultiNEST \cite{feroz2009multinest}, NS-SMC \cite{salomone2018unbiased}, dynamic nested sampling \cite{higson2019dynamic} or others), but a preprocessing step after which such a modified nested sampling should be carried out on the non-plateau component. In more detail, we propose the following approach:
\begin{enumerate}
    \item generate $m$ samples $x_i, i = 1,...,M$ from measure $\mu$ and sort them:
           \begin{itemize}
               \item $x_i \in \Omega_r$: sample becomes part of active set for nested sampling (in total $n_r$ active samples)
               \item $x_i \in \Omega_{h,c}$: sample used for estimation of $\Delta_h$
            \end{itemize}
            
    \item calculate evidence
        \begin{enumerate}
            \item estimate $\Delta_h$ via relative frequencies: $\Delta_h = \frac{n_h}{m}$ with $n_h$ number of samples $x_i \in \Omega_{h,c}$
            
            \item compute $Z_c = \sum_{h} \L_{h,c} \cdot \Delta_h$
            
           \item perform nested sampling on $\Omega_r$ considering the maximum accessible prior mass, i.e. $\xi_0 = 1-\sum_h \Delta_h$ (instead of $\xi_0 = 1$) and obtain $Z_r$ as depicted in Figure~\ref{fig:decompositionEvidence}

            \item calculate evidence $Z \approx Z_r + Z_c$ according to \eqref{eq:robustImplementation_calculationEvidence}
            \end{enumerate}
\end{enumerate}

This preprocessing method allows also for the approximation of mean and variance of the posterior measure $\mu_\text{post}$  and the generation of weighted samples from it.

\subsection{Numerical experiments: Randomization strategy versus splitting approach}

 In the case of likelihood plateaus with large enough prior measure, it is ineffective to run nested sampling without any consideration for the plateau (i.e.\ without randomization and the splitting method). The reason for that is that a fundamental assumption of nested sampling, the approximation $X_k \approx \e^{-k/n}$, is not valid. This can be seen from the cumulative distribution function for $\mu \circ \Phi^{-1}$ which differs from the cumulative distribution function of a uniform distribution. Furthermore, in a setting where the plateau is at the maximum of the likelihood, an additional difficulty arises: As soon as all active samples are within the plateau, there is no ``better'' sample in terms of the likelihood value. Thus, nested sampling is forced to shut down before its convergence and this results in insufficient numerical approximation of the evidence.


To overcome this type of situation, it has been suggested to use a randomization strategy \cite{skilling2006nested} meaning that every evaluation of the likelihood during nested sampling is independently perturbed by a random variable $\xi$ resulting in
\begin{align}
    \L_n^+ (x) = \L(x) + \xi_n,
\end{align}
where $\L_n^+$ is the $n$-th evaluation of the likelihood and $\xi_n$ is i.i.d. perturbative noise, for example a Gaussian sample or a uniformly distributed random variable.

\begin{figure}
    \centering
    \includegraphics[width=0.5\textwidth]{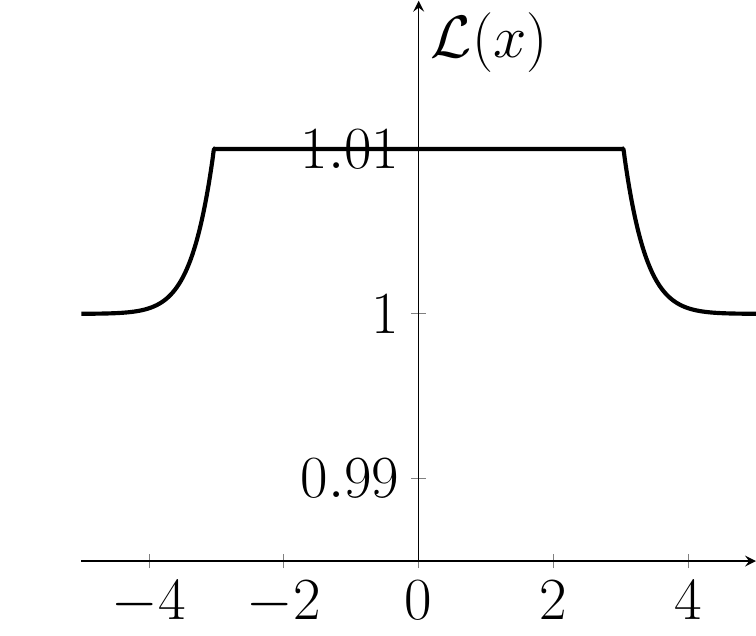}
    \caption{Likelihood function $\mathcal{L}(x)$ along $x_1$-axis.}
    \label{fig:likelihood}
\end{figure}

We compare our proposed splitting approach to this randomization strategy with $\xi$ being uniformly distributed. We consider the parameter space $\R^5$ with the prior being chosen as $\mu = \mathcal N(0_5, 2^2\cdot E_5))$ where $0_5$ is the zero vector in $\R^5$ and $E_5$ is the identity matrix in $\R^{5\times 5}$. As likelihood function, we use
\begin{align}
    \L(x) = \min\{1 + \e^{-\frac{\lVert x \rVert^2}{2}}, \L_\mathrm{plateau} \} \label{eq:refSettingNumericalExperiments}
\end{align}
with $\L_\mathrm{plateau} = 1.01$. This means that $\L$ is a capped normal Gaussian without normalization (see Figure \ref{fig:likelihood}). The plateau covers approximately 19.4\% of the prior mass, i.e.\ $X(\L_\mathrm{plateau}) \approx $. We additionally increase the likelihood level by $+1$ so we can guarantee that the perturbed likelihood evaluation $\L_n^+$ remains positive even under additive perturbation by $\xi$. This is also the reason why we choose a uniformly distributed (as compared to a normally distributed and thus unbounded) perturbation, i.e.\,$\xi \sim \operatorname{Unif}[-a,a]$ with $a = \left\{10^{-3}, 10^{-5}, 10^{-7}\right\}$. The random variations in the likelihood correspond, thus, to 10\%, 0.1\% and 0.001\% of the total range of the likelihood, which is the interval $[1.00, 1.01]$. The analytical evidence is $Z_a \approx 1.002694$. 

Since in many applications the evaluation of the likelihood function is extremely costly, we use a budget of likelihood evaluations in the numerical experiments. We relate the number of likelihood evaluations to the number of total active samples $n$ and we define a low budget regime ($10\cdot n$ likelihood evaluations), and a high budget regime ($40 \cdot n$ likelihood evaluations). At the same time, we compare the quality of the approximation of the evidence via both methods (perturbation and splitting) for three different choices of active samples ($n = 50, 100, 500$). 

\begin{figure}
    \centering
    \includegraphics[width=0.7\textwidth]{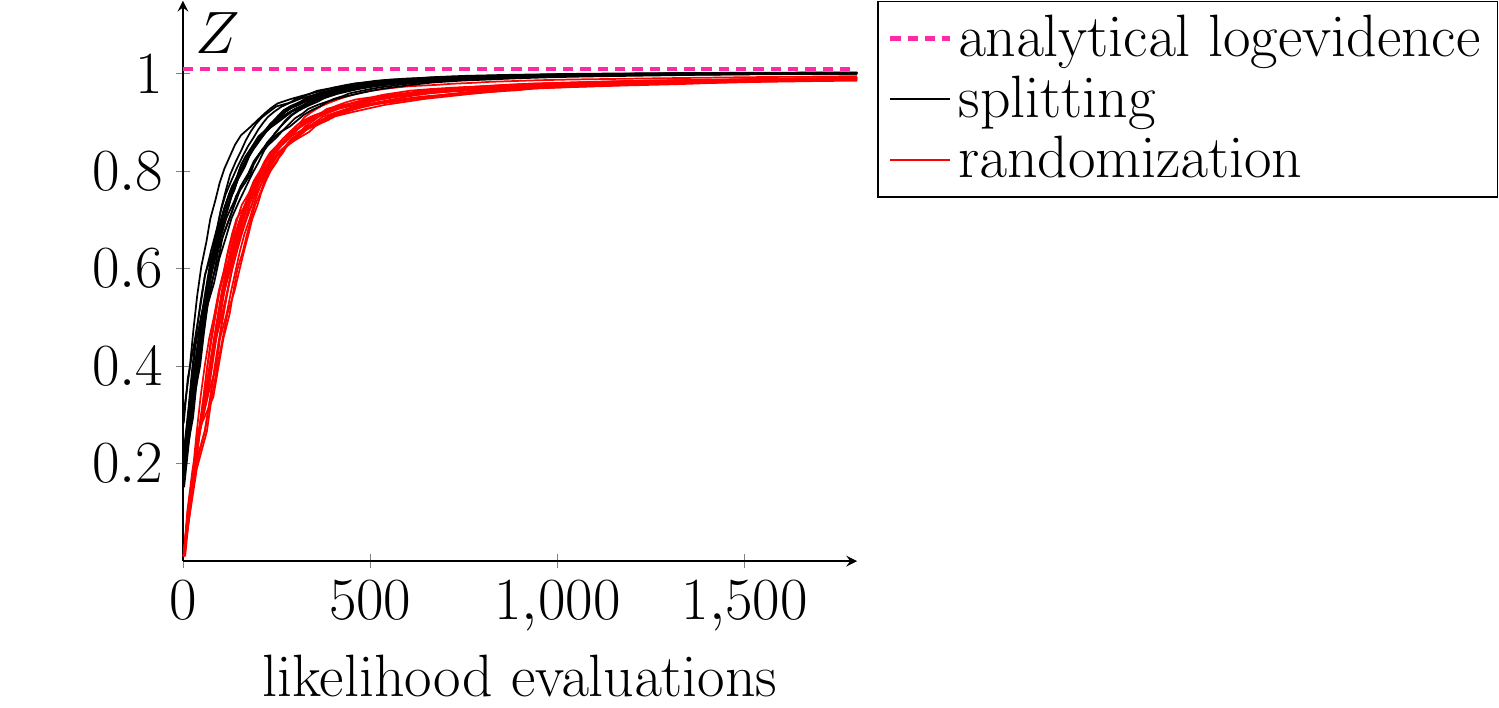}
    \caption{Evolution of the evidence over likelihood evaluations for the 5D setting with 100 active samples ($n = 100$) and the likelihood randomization $\operatorname{Unif}[-10^{-5}, 10^{-5}]$. Shown are 15 representatives for each approach. The evidence obtained from the splitting method converges faster to the analytical evidence in terms of likelihood evaluations than the evidence obtained from the randomization method. }
    \label{fig:SkillingVSsplitting_EvidenceEvolution}
\end{figure}

Most importantly, the numerical computation of the evidence converges to the analytical solution with both methods, i.e. via the splitting and the randomization method. Representative curves illustrating the evolution of the evidence per likelihood evaluations are given in Figure~\ref{fig:SkillingVSsplitting_EvidenceEvolution} for 100 active samples.

\begin{figure}
    \centering
    \includegraphics[width=\textwidth]{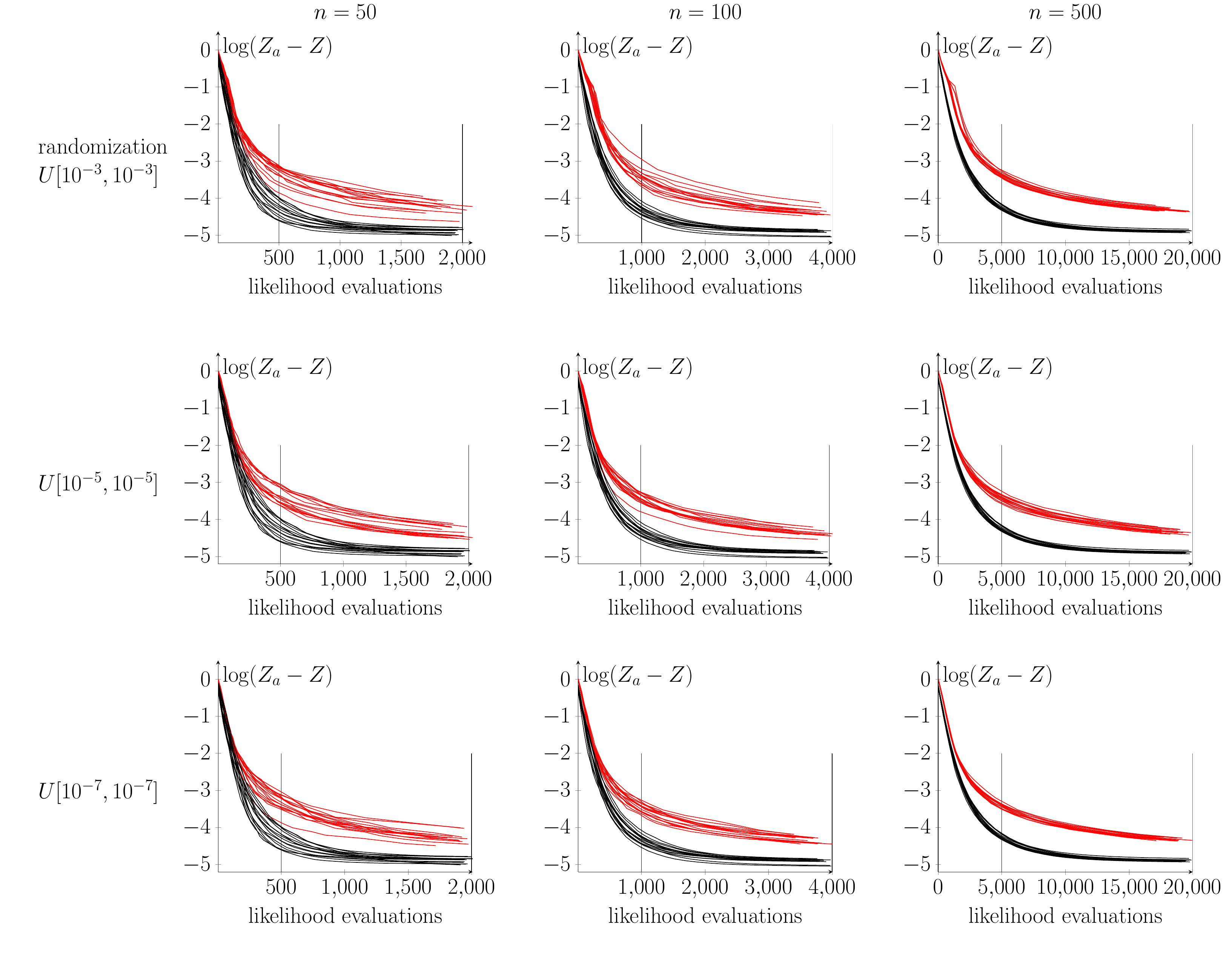}
    \caption{Logarithmic error of evidence estimation $\log(Z_a - Z)$ over likelihood evaluations for 15 representative runs for each parameter set $n = \left\{50, 100, 500\right\}$ and randomization $a = \left\{ 10^{-3}, 10^{-5}, 10^{-7}\right\}$. The black curves obtained from splitting method do not vary in row but in column whereas the red curves of the randomization strategy changes in dependes of $a$.}
    \label{fig:SkillingVSsplitting_logError}
\end{figure}

The fact that the contribution of the likelihood plateaus to the evidence is approximated very quickly before starting nested sampling on the remainder of the domain results in the splitting method starting with an initial lead of (in this example) approximately $0.2$, which is exactly the numerical value of the evidence contributed by the likelihood plateau. We note that this lead is not bought with additional computational overhead as vanilla nested sampling requires initial evaluation of the likelihood in all initial samples, too.
But it is not just this head start leading to a relative improvement of the splitting method over perturbed nested sampling, as can be seen from the logarithmic error plots in figure \ref{fig:SkillingVSsplitting_logError}, which shows a quicker rate of approximation as well.

The splitting method converges faster in terms of likelihood evaluations for all $n$ and all magnitudes of randomization we explored. This can be seen from the diagrams in Figure~\ref{fig:SkillingVSsplitting_logError} showing the logarithmic approximation error of the evidence given by $\log(Z_a - Z)$. Here, the black curves representing the splitting method reach a specific error within a much lower number of likelihood evaluations. The spread of the curves decreases with increasing $n$ as to be expected for a large number of active samples. 

\begin{figure}
    \centering
    \includegraphics[width=\textwidth]{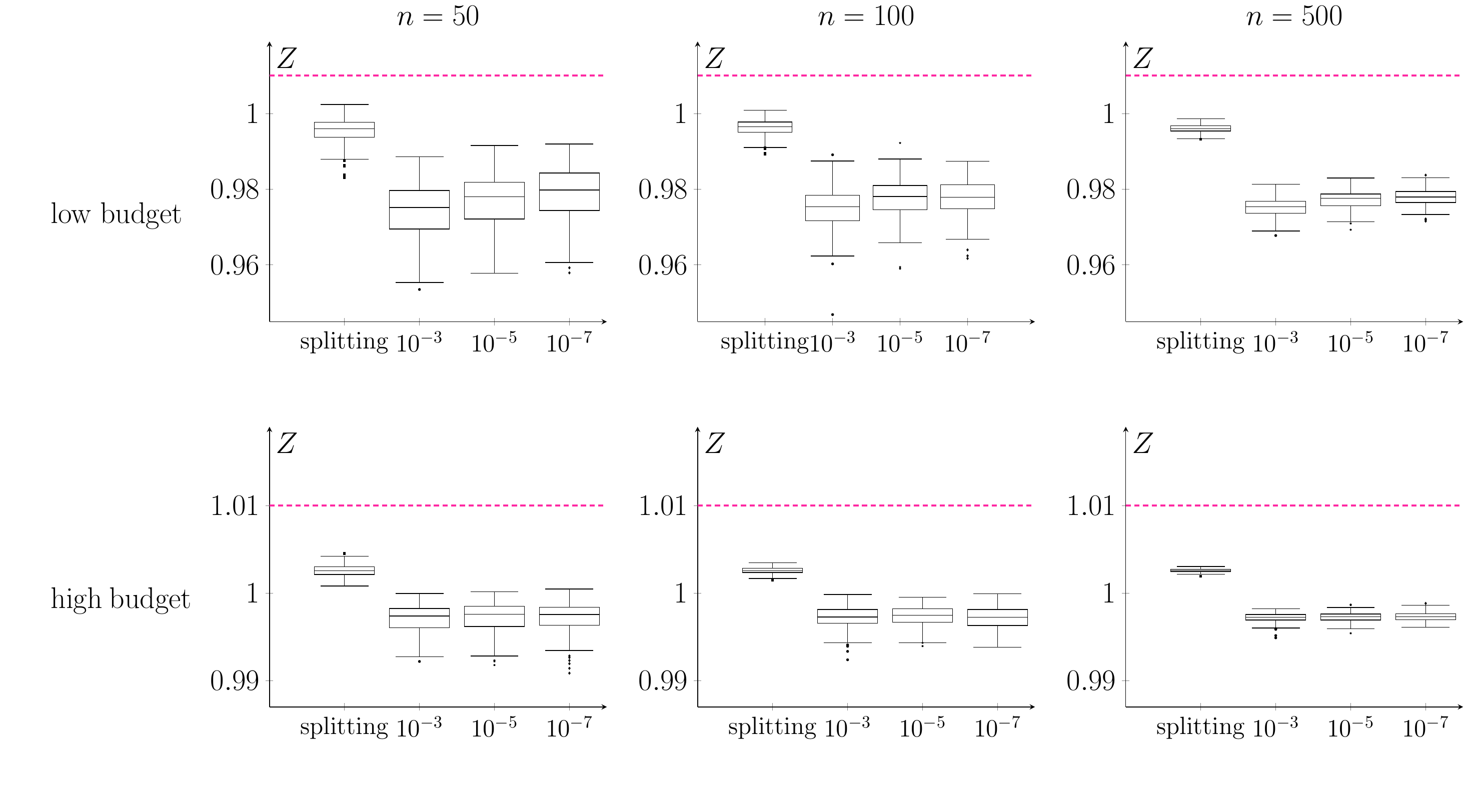}
    \caption{Boxplot of estimated evidences of 200 runs using the splitting method and the randomization strategy with $a = \left\{10^{-3}, 10^{-5}, 10^{-7}\right\}$. Using a high budget of likelihood evaluations ($40\cdot n$) leads to a better approximation of the evidence than a low budget ($10 \cdot n$). The analytical evidence $Z_a = 1.002694$ is represented by the dashed line.}
    \label{fig:SkillingVSsplitting_Boxplot}
\end{figure}

The spread of evidence computation between independent runs of both methods is further explored in Figure~\ref{fig:SkillingVSsplitting_Boxplot}. Here, evidence approximations of 200 independent runs for each parameter set are visualized. It is apparent that the choice of perturbation magnitude $a$ has some influence on the quality of evidence approximation: the lower $a$, the better the approximation, although even the smallest perturbations do not reach the quality of approximation reached by the splitting method in the same amount of likelihood evaluations. This is true for both the low-budget and high-budget regime.

In total, the numerical experiments show that the proposed pre-processing step requires a lower number of likelihood evaluations for a specific approximation error compared to the randomization method. Of course, this effect scales with the magnitude with which the plateau contributes to the evidence: the higher the plateau contribution, the more effective the splitting approach. And, of course, in settings where the likelihood plateaus carry only an insignificant amount of prior mass, the splitting approach will only slightly improve the overall quality of approximation.

The remainder of the manuscript is devoted to the proof of theorem \ref{thm:main-theorem} and can be safely skipped by readers mainly interested in the practical application of nested sampling.
\begin{appendix}
\section{The pushforward measure $\Phi_\#\mu$}
\label{sec:uniform}
In this section we will talk about the equation
\begin{equation}\label{eq:push-forward-true}
    \int_{0}^1 \tilde{\mathcal{L}}(r) \diff (\Phi_\#\mu)(r)= \int_0^1\tilde{\mathcal{L}}(X) \diff X.\end{equation}

In the simplest case, this would be true if $\Phi_\#\mu$ is the uniform measure on $[0,1]$, because then $\frac{\diff \Phi_\#\mu(r)}{\diff r} = 1$. This is equivalent to saying that i.i.d $\mu$-samples $\{x_i\}_i\in\Omega$ become uniform samples $\{\Phi(x_i)\}_i$ in $[0,1]$ via the mapping $\Phi$. 

We will show that $\Phi_\#\mu = \operatorname{Unif}[0,1]$ if and only if there are no plateaus of $\L$ which have positive $\mu$-measure (i.e.\ $\L_\#\mu$ has not atoms). 

This means that in the presence of plateaus in $\L$, equation~\eqref{eq:push-forward-true} is more difficult to prove (because integration with respect to $\Phi_\#\mu$ is not standard Lebesgue integration), but we will nevertheless be able to show correctness of this equation.

The cumulative distribution function of $\Phi_\#\mu$ is
	\begin{align} \label{eq:equivformuniform}
		\Phi_\#\mu ([0,\alpha]) &= \mu \left( \left\lbrace x \in \Omega: \Phi(x) \in [0,\alpha] \right \rbrace  \right)\\ \notag
		&= \mu \left( \left\lbrace x \in \Omega : X(\mathcal{L}(x)) \in [0,\alpha] \right\rbrace \right)\\ \notag
		&= \mu \left( \left\lbrace x \in \Omega : \mu \left( \left\lbrace z \in \Omega : \mathcal{L}(z) > \mathcal{L}(x) \right\rbrace \right) \le \alpha \right \rbrace\right)
	\end{align}
This means that  $\Phi_\#\mu = \operatorname{Unif}[0,1]$ if and only if
\begin{equation}\label{eq:mainequationuniform}
    \mu \left( \left\lbrace x \in \Omega : \mu \left( \left\lbrace z \in \Omega : \mathcal{L}(z) > \mathcal{L}(x) \right\rbrace \right) \le \alpha \right \rbrace\right) = \alpha.
\end{equation}

In this section we derive conditions for the validity of \eqref{eq:mainequationuniform}. Key statement is the following lemma.

\begin{lemma}
\label{lem:main}
Let $\mu$ be a measure on a measurable space $(\Omega, \Acal)$, $\Omega \subseteq \mathbb{R}^d, d\in \mathbb{N}$, and $\L : (\Omega, \Acal) \to (\R, \mathcal B)$ a measurable mapping which is bounded from above and below, i.e.\ $-\infty < \inf \L \leq \sup \L < \infty$. We consider the push-forward of $\mu$ via $\L$, denoted by $\L_\#\mu$ (which now is a measure on $(\R, \mathcal B)$). We assume that we know all positions $r_i\in \R$ at which the discrete part of $\L_\#\mu$ has a non-vanishing contribution, in particular we assume that we know the $(r_i)_{i=1}^N, (\alpha_i)_{i=1}^N,$ and $(\Delta_i)_{i=1}^N$ such that
\begin{align*}
    \Lmu(\{r_i\}) &= \mu(\{z\in \Omega: \L(z) = r_i\}) = \Delta_i\\
    \Lmu((r_i, \infty)) &=  \mu(\{z\in \Omega: \L(z) > r_i\}) = \alpha_i.
\end{align*}
We assume that there are only finitely many (i.e.\ $N$) of such atoms and that $\Lmu$ is non-singular everywhere else. Then for any $\alpha \in [0, 1]$,
\begin{align*} &\mu\left(\left\{x\in \Omega:\mu\left(\left\{z\in \Omega: \L(z) > \L(x)\right\}\right)\leq \alpha\right\}\right)\\
&\qquad= \begin{cases} \alpha_i + \Delta_i &\text{for } \alpha \in [\alpha_i, \alpha_i + \Delta_i)\\ \alpha &\text{else.}
\end{cases}\end{align*}
For a visualization of this and the quantities $\alpha_i, \Delta_i$, see Fig.~\ref{fig:L-and-psiralpha}.


\begin{figure}
    \centering
    \includegraphics[width=1\textwidth]{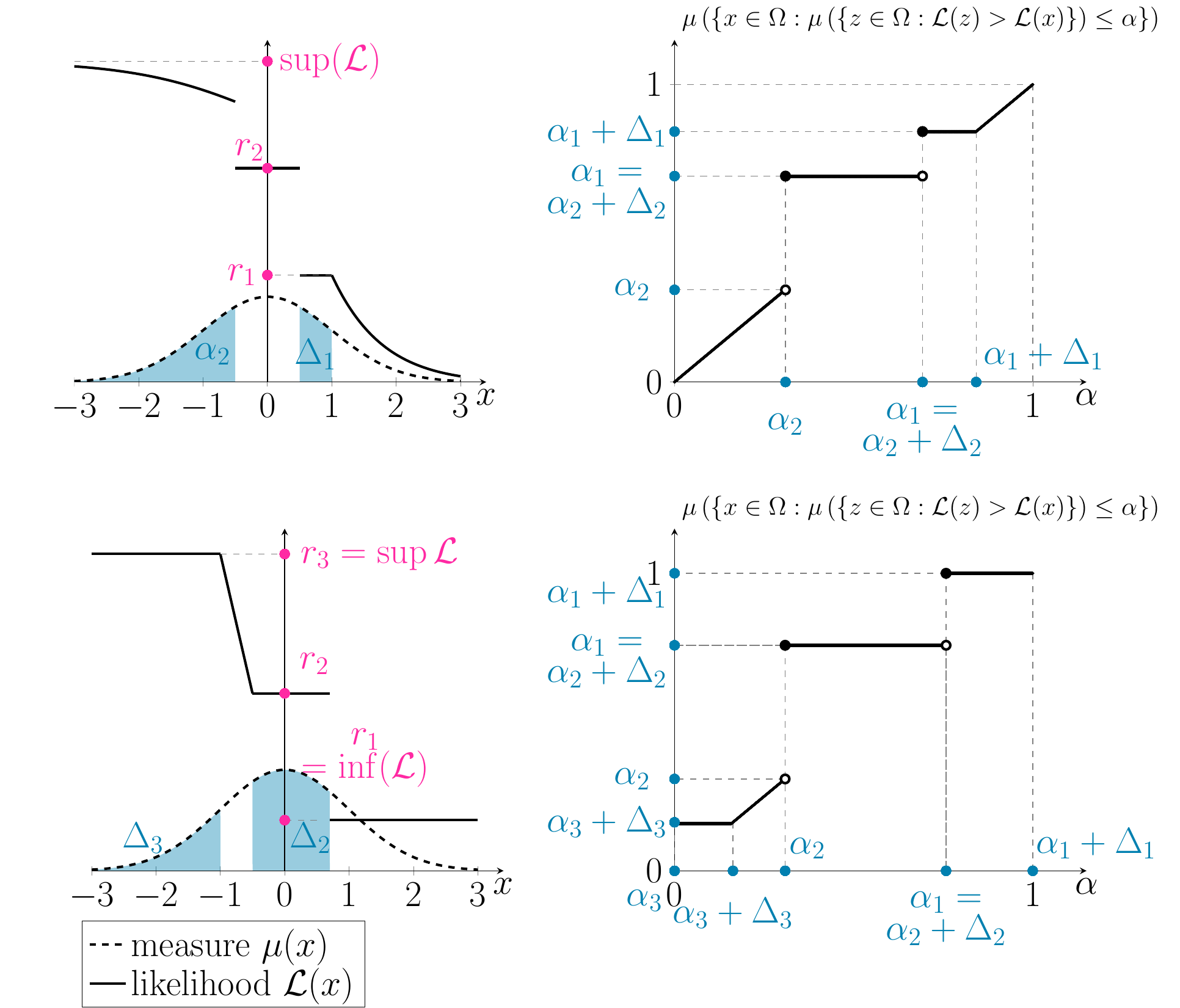}
    \caption{Row-wise examples of measure $\mu$ and likelihood function in the context of Lemma~\ref{lem:main}. Apparently, $\mu\left(\left\{x\in \Omega:\mu\left(\left\{z\in \Omega: \L(z) > \L(x)\right\}\right)\leq \alpha\right\}\right)$ is not uniform e.g.\ if $\L$ has plateaus on the support of $\mu$, neither in case 0 (row 1) nor in case AB (row 2).}
    \label{fig:L-and-psiralpha}
\end{figure}

\begin{remark}

It was suggested to us that (in order to reduce the complexity of the proof) we use the existing machinery of generalized inverses of the survival function (or quantile function). Employing this theory would necessitate some heavy translational work which would invariable impede ease of reading. For example, the function for which we will need to construct a generalized inverse will be non-increasing which goes against the convention of generalized inverses which are constructed for non-decreasing functions. After having done this compatibility conversion, we could then ``outsource'' Step 2 of the following proof by citing respective results about generalized inverses (and translating back to our anti-parallel setting). Unfortunately this would ultimately only reduce the length of the proof by a small fraction, but we would need to sacrifice expository completeness in the process. The bulk of the (admittedly technical) proof will be paying attention to various edge cases and other matters not connected to the theory of generalized inverses. For this reason we have decided to present the proof without excursion into this field.
\end{remark}

\begin{proof}[Proof of Lemma~\ref{lem:main}]

We sketch the proof formally first: Key ingredients will be two functions which, concatenated, yield the result. First, we will define $\psi(r) = \mu(\{z\in \Omega: \L(z)\geq r\}$ and then we define $\alpha\mapsto r_\alpha = \inf \{r \in \mathbb{R}: \psi(r)\leq \alpha\}$ (a kind of generalized inverse for $\psi$). Note that $\psi$ is similar to the function $X$ defined in Section~\ref{sec:mainParadigm}, but with a $\geq$ instead of $>$, which we will need. Furthermore, $r_\alpha = \tilde \L(\alpha)$, although this is not obvious at this point and will be proven in step~\ref{step:definition-ralpha}.\ref{sstep:tildeL}. We will be able to prove that 
\[\psi(r_\alpha) = \begin{cases} \alpha_i + \Delta_i &\text{for } \alpha \in [\alpha_i, \alpha_i + \Delta_i)\\ \alpha &\text{else.}
\end{cases}\]
and also that $\mu(\{x\in \Omega: \mu(\{z\in \Omega: \L(z) > \L(x)\}) \leq \alpha\}) = \mu(\{x\in \Omega: \L(x) \geq r_\alpha\}) = \psi(r_\alpha)$, which, in combination, yield the result.

Due to a few technicalities, there are some pathological cases that need to be taken care of individually.

We show the result in a series of steps. First, we will need to fix notation for the actual proof.

We start by ordering $r_1 < r_2 < \cdots$. Note that two $r_i, r_j$, $i\neq j$ cannot be identical or we could just drop one of them. With this ordering, we obtain $\alpha_1 + \Delta_1 > \alpha_1 \geq \alpha_2 + \Delta_2 > \alpha_2 \geq \cdots$, i.e.
\begin{align*}
    \alpha_{i} + \Delta_i &> \alpha_i\\
    \alpha_i &\geq \alpha_{i+1} + \Delta_{i+1}
\end{align*}
Indeed, $\alpha_k + \Delta_k > \alpha_k$ because $\Delta_k > 0$ and $\alpha_k = \Lmu((r_k, \infty)) \geq \Lmu([r_{k+1}, \infty)) = \alpha_{k+1} + \Delta_{k+1}$ because $r_{k+1} > r_k$. Here, it is possible to have equality, i.e.\ a situation where $\Lmu((r_k,\infty)) = \Lmu([r_{k+1},\infty))$ or equivalently $\alpha_k = \alpha_{k+1}+\Delta_{k+1}$.
 
 From this point on, we will without loss of generality always consider this ordering of the $r_i$ and $\alpha_i$.

There are two special cases to consider depending on whether the infimum and the supremum of $\L$ constitute plateaus with positive measure themselves: Special Case A will be when $r_1 = \inf \L$ which holds if and only if $\alpha_1 + \Delta_1 = 1$. Special Case B is when $r_N = \sup \L$ (i.e.\ if and only if $\alpha_N = 0$). Any of those two cases can be present (similar to blood types), i.e.\ we will always have four cases to keep in mind: A, B, AB, and 0 (i.e.\ neither A nor B).

Before we can continue, we need to define specific decompositions for the intervals $[\inf\L, \sup \L]$ and $[0,1]$ (the range of possible probabilities). 

We set \[I_i = (r_i, r_{i+1}], \quad, i=1,\ldots,N-1.\]

The leftmost and rightmost intervals need to be chosen according to the special case $\in$ (A, B, AB, 0) we are in: 
\[
I_0 = \begin{cases} \{\inf \L\} ,&\text{ if } r_1 = \inf \L, \text{ i.e.\ in cases A, AB}\\
[\inf \L, r_1],&\text{ if } r_1 > \inf \L, \text{ i.e.\ in cases 0, B}
\end{cases}
\]

In the same way, we need to be careful with the other boundary. If $r_N = \sup \L$, then $I_{N-1} = (r_{N-1}, r_N]$ is the last interval. If $r_N < \sup \L$, we need to additionally define $I_N = (r_N, \sup \L]$. In order to simplify notation, we write
\[
I_N = \begin{cases} \emptyset ,&\text{ if } r_N = \sup \L, \text{ i.e.\ in cases B, AB}\\
(r_N, \sup \L],&\text{ if } r_N < \sup \L, \text{ i.e.\ in cases 0, A}
\end{cases}
\]

At any rate, $[\inf \L,\sup \L] = \bigcup_{i=0}^N I_i$ is a disjoint decomposition (see Fig.~\ref{fig:decomposition-interval-infL-supL}).

Similarly we define a decomposition of the interval $[0,1]$: We set \[A_i = [\alpha_i, \alpha_i + \Delta_i).\] for $i=1,\ldots, N$. The ``missing parts'' are defined via  \[J_i=\begin{cases}[\alpha_{i+1}+\Delta_{i+1}, \alpha_i), &\text{ if } \alpha_{i+1}+\Delta_{i+1} < \alpha_i \\ \{\alpha_i\},&\text{ if } \alpha_{i+1}+\Delta_{i+1}= \alpha_i \end{cases}\] for $i=1,\ldots, N-1$.  

As before, we need to pay attention to details concerning the ``boundary conditions'':

\[
J_0 = 
\begin{cases}
\{1\}, &\text{ if } r_1 = \inf \L, \text{ i.e.\ in cases A, AB}\\
[\alpha_1+\Delta_1, 1], &\text{ if } r_1 > \inf \L, \text{ i.e.\ in cases 0, B}
\end{cases}
\]
and 
\[
J_N = \begin{cases} \emptyset ,&\text{ if } r_N = \sup \L, \text{ i.e.\ in cases B, AB}\\
[0, \alpha_N),&\text{ if } r_N < \sup \L, \text{ i.e.\ in cases 0, A}
\end{cases}
\]
A visualization of the intervals is shown in Fig.~\ref{fig:decomposition-interval-alpha}.

\begin{step}\label{step:intervalDefinition}
The intervals $I_i$ yield a complete and disjoint decomposition\footnote{We use $\cup$ for union of sets in general and $\sqcup$ explicitly for disjoint unions of sets, i.e.\ $A\sqcup B$ contains no points which are in both $A$ and $B$.}
\[[\inf\L, \sup\L] = \bigsqcup_{i=0}^N I_i.\]
The intervals $A_i = [\alpha_i, \alpha_i + \Delta_i)$ are disjoint from each other and are ``proper'' intervals, i.e.\ $\alpha_i + \Delta_i > \alpha_i$. The $A_i$ and the $J_i$ together form a decomposition:
\[[0,1] = J_0 ~\sqcup~ \bigsqcup_{i=1}^N (A_i \cup J_i).\]
Note that this decomposition is not completely disjoint: $A_i\cap J_i = \{\alpha_i\} = J_i$ if and only if $\alpha_{i+1}+\Delta_{i+1} = \alpha_i$. Apart from that, all sets involved are disjoint from each other, in particular $(A_i \cup J_i) \cap (A_j \cup J_j) = \emptyset$ for $i\neq j$.
\end{step}

\begin{proof}[Proof of step~\ref{step:intervalDefinition}]
 Proven in the text above.
\end{proof}

\begin{figure}
    \begin{subfigure}{\textwidth}
        \centering
        \includegraphics[width=0.8\textwidth]{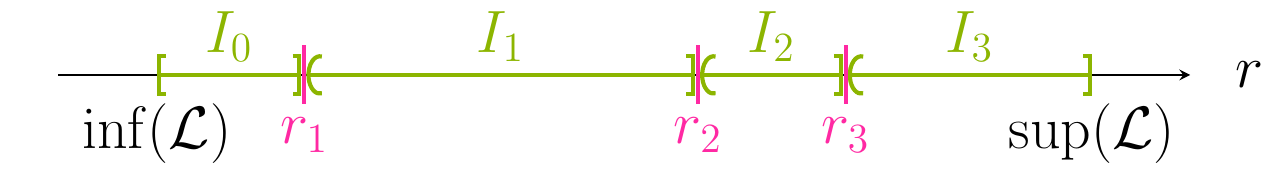}
        \caption{case 0.} \label{fig:decomposition-interval-infL-supL-case0}
    \end{subfigure}
    \begin{subfigure}{\textwidth}
        \centering
        \includegraphics[width=0.8\textwidth]{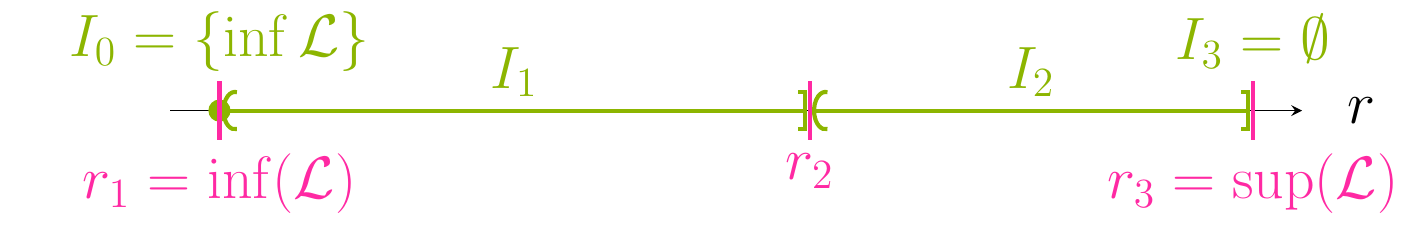}
        \caption{case AB.} \label{fig:decomposition-interval-infL-supL-caseAB}
    \end{subfigure}
    \caption{Decomposition of interval $[\inf \L, \sup \L]$ into $I_i$.}
    \label{fig:decomposition-interval-infL-supL}
\end{figure}
\begin{figure}
    \begin{subfigure}{\textwidth}
        \centering
        \includegraphics[width=0.8\textwidth]{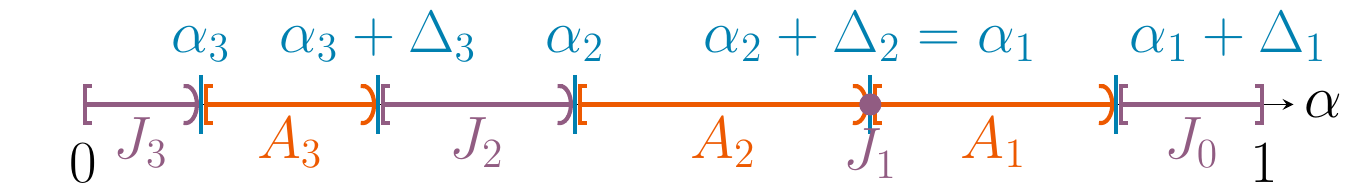}
        \caption{case 0.} \label{fig:decomposition-interval-alpha-case0}
    \end{subfigure}
    \begin{subfigure}{\textwidth}
        \centering
        \includegraphics[width=0.8\textwidth]{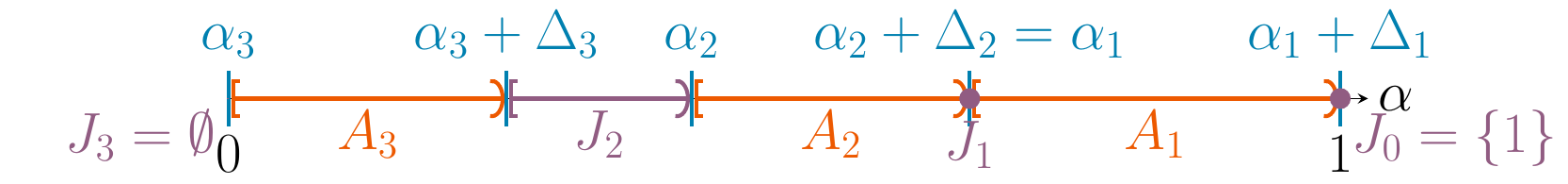}
        \caption{case AB.} \label{fig:decomposition-interval-alpha-caseAB}
    \end{subfigure}
        \caption{Decomposition of interval $[0,1]$ into $A_i$ and $J_i$.}
        \label{fig:decomposition-interval-alpha}
\end{figure}
   
\begin{step} \label{step:mappingPsi}
We consider first the map $\psi: r\mapsto \mu(\{z\in \Omega: \L(z) \geq r\}) = \Lmu([r,\infty))$. This mapping has the following properties:
    \begin{enumerate}
    \item \label{sstep:psi-nonincreasing}
        $\psi$ is non-increasing.
    \item  \label{sstep:psi-leftcontinuous}
      $\psi$ is left-continuous at $r=r_i$, with $\lim_{r\nearrow r_i}\psi(r) = \psi(r_i) = \alpha_i+\Delta_i$.
    \item \label{sstep:psi-nonrightcontinuous}
        $\psi$ is non-right-continuous at $r=r_i$, because we have $\lim_{r\searrow r_i} \psi(r) = \alpha_i < \alpha_i+\Delta_i = \psi(r_i)$.
    \item \label{sstep:psi-i-surjective}
        Any restricted mapping $\psi_i = \psi|_{I_i}$ is a continuous, well-defined, surjective, and non-increasing map $\psi_i: I_i \to J_i$ . More precisely we can say the following:
        \begin{itemize}
          \item If $i=0$, then either the mapping is trivial (in cases A or AB) in the sense that $\psi_0: \{\inf\L\} \to \{1\}$ or (in cases 0 or B): $\psi_0 : I_0 = [\inf \L, r_1] \to [\alpha_1+\Delta_1, 1]$ is a continuous, non-increasing and onto mapping with $\psi(r_1) = \alpha_1+ \Delta_1$ and $\psi(\inf \L) = 1$.
         \item  If $i=N$, then either the mapping is trivial (in cases B or AB) in the sense that  $\psi_N: \emptyset \to \emptyset$ or (in cases 0, A): $\psi_N:(r_N, \sup\L] \to [0, \alpha_N)$ is a continuous, non-increasing and onto mapping with $\lim_{r\searrow r_N}\psi(r) = \alpha_N$ and $\psi(\sup \L) = 0$.
            \item If $i=1,\ldots, N-1$, then there are two possibilities:
            \begin{itemize}
       
                \item $\alpha_{i+1}+\Delta_{i+1} = \alpha_i$ and $\psi_i \equiv \alpha_i$, i.e.\ $\psi_i: I_i \to J_i$ collapses $I_i$ to one point, or
                \item $\alpha_{i+1}+\Delta_{i+1} < \alpha_i$. Then $\psi_i$ is continuous, non-increasing and surjective as a function $I_i \to J_i$.
        
                If $i = 1,\ldots, N-1$, then: $\lim_{r\searrow r_i}\psi(r) = \alpha_i$ and $\psi(r_{i+1}) = \alpha_{i+1}+\Delta_{i+1}$. 
            \end{itemize}
        \end{itemize}
        These properties can also be seen in Fig.~\ref{fig:psiOfr-intervalConnectionPsi-case0} and Fig.~\ref{fig:psiOfr-intervalConnectionPsi-caseAB} from our example at the beginning (\ref{fig:L-and-psiralpha}).
        
    \item \label{sstep:phi}
    (For later purposes): We define a slight variation on $\psi$ given by $\phi: r\mapsto \mu(\{z\in \Omega: \L(z) > r\}) = \Lmu((r,\infty))$. Then $\phi$ is the right-continuous version of $\psi$, i.e.\ $\phi(r) = \psi(r)$ for any $r\neq r_i$ and $\phi(r_i) = \alpha_i$. Note that $\phi(r) = X(r)$.
    \end{enumerate}

\begin{figure}
\centering
    \begin{subfigure}{0.7\textwidth}
        \centering
        \includegraphics[width=1\textwidth]{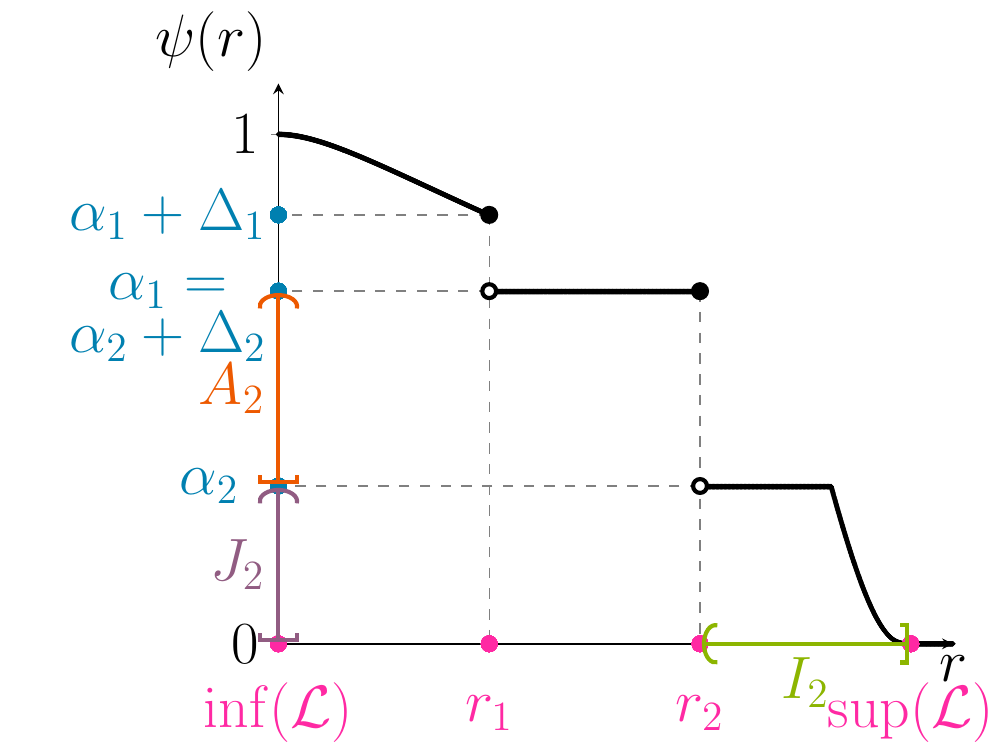}
        \caption{Mapping $\psi(r)$.} \label{fig:psiOfr-case0}
    \end{subfigure}\\
    
        \begin{subfigure}{1\textwidth}
        \centering
        \includegraphics[width=0.9\textwidth]{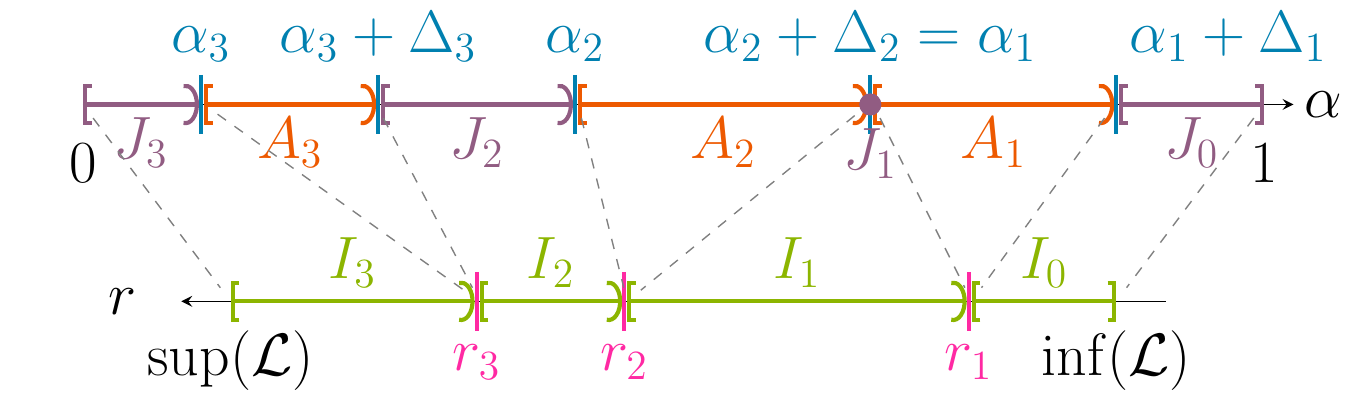}
        \caption{Interval mapping via $\psi(r)$.} \label{fig:intervalConnectionPsi-case0}
    \end{subfigure}
    
    \caption{Case 0: $\psi(r)$ and interval mapping from $[\sup\L,\inf\L]$ to $[0,1]$ via  $\psi(r)$ for the likelihood and the measure shown in figure \ref{fig:L-and-psiralpha} (1\textsuperscript{st} row).} \label{fig:psiOfr-intervalConnectionPsi-case0}
\end{figure}

\begin{figure}
\centering
\begin{subfigure}{0.7\textwidth}
        \centering
        \includegraphics[width=1\textwidth]{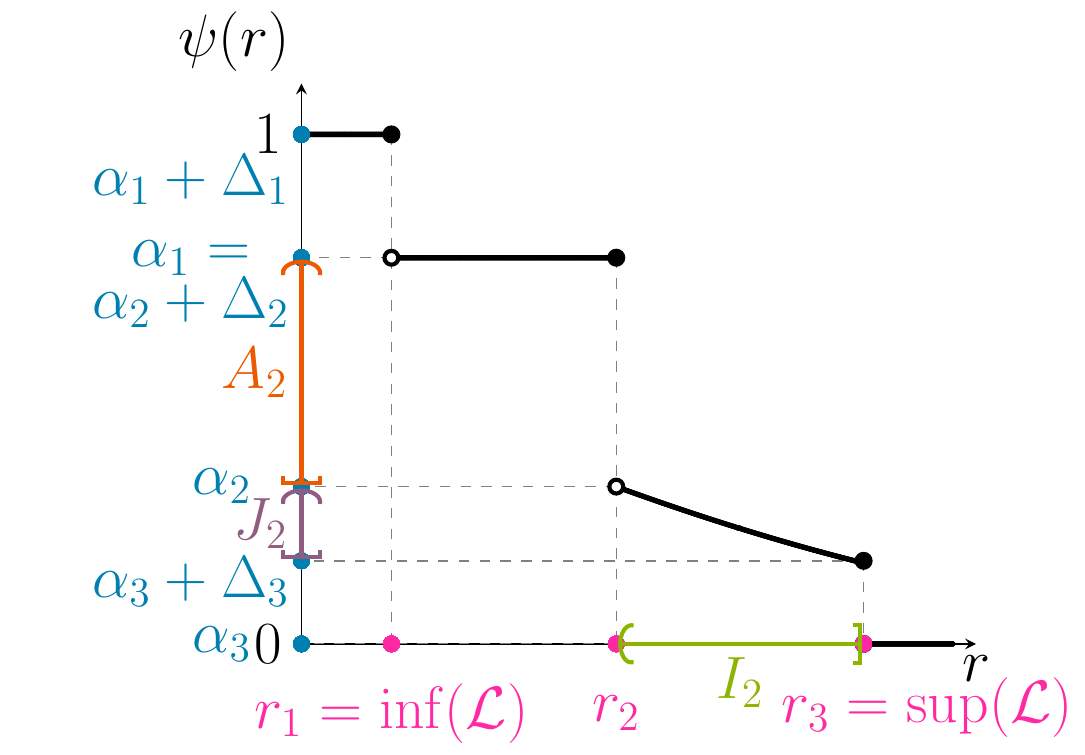}
        \caption{Mapping $\psi(r)$.} \label{fig:psiOfr-caseAB}
    \end{subfigure}
\\
    \begin{subfigure}{1\textwidth}
        \centering
        \includegraphics[width=0.9\textwidth]{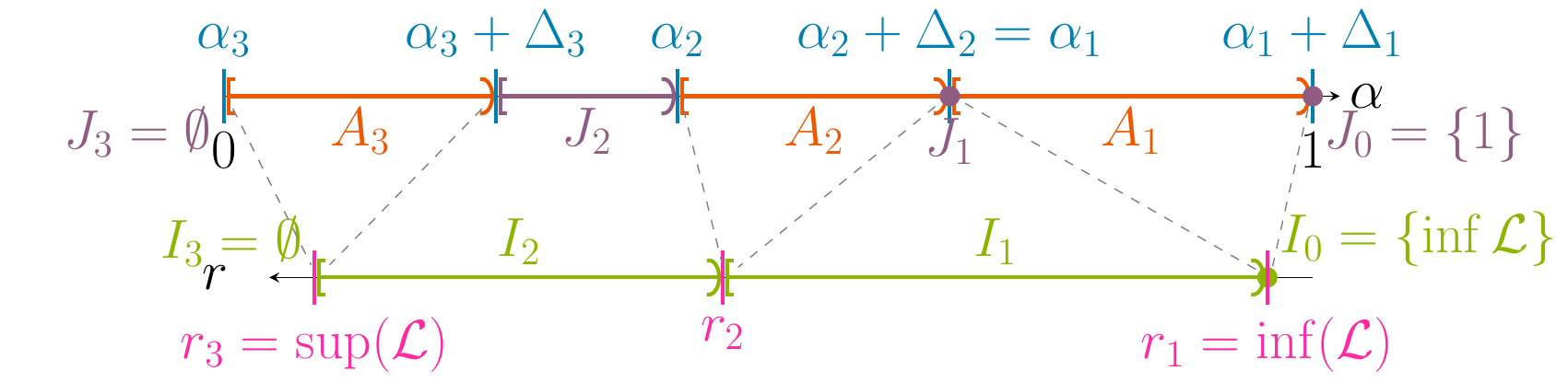}
        \caption{Interval mapping via $\psi(r)$.} \label{fig:intervalConnectionPsi-caseAB}
    \end{subfigure}
    
    \caption{Case AB: $\psi(r)$ and interval mapping from $[\sup\L,\inf\L]$ to $[0,1]$ via $\psi(r)$ for the likelihood and the measure shown in Fig.~\ref{fig:L-and-psiralpha} (2\textsuperscript{nd} row).}\label{fig:psiOfr-intervalConnectionPsi-caseAB}
\end{figure}
\end{step} 

\begin{proof}[Proof of step~\ref{step:mappingPsi}]
 Ad \ref{sstep:psi-nonincreasing}: Define $B_r = \{z\in \Omega: \L(z) \geq r\}$ The monotonicity of $\psi$ follows from the inclusion of the sets $B_s \subseteq B_r$ for $r \leq  s$. 

 Ad \ref{sstep:psi-leftcontinuous}: Let $(r^{(n)})_{n\in\N}$ be any sequence monotonously increasing with limit $r_i$ for some $i$, i.e.\ $r^{(n)}\nearrow r_i$. We need to show that $\lim_{n\to\infty} \psi(r^{(n)}) = \psi(r_i) = \alpha_i+\Delta_i$. Define $E_n = \{z\in \Omega: \L(z) \geq r^{(n)}\}$, and thus $E_n \supset E_{n+1}\supset \cdots$ and $\bigcap_n E_n = \{z\in \Omega: \L(z) \geq r_i\}$. The measure $\mu$ is continuous from above, which implies the identity marked with $\star$ in the identity $\lim_{n\to\infty}\psi(r^{(n)}) = \lim_{n\to\infty}\mu(E_n)  \stackrel{\star}{=} \mu(\bigcap_n E_n) = \mu(\{z\in \Omega: \L(z) \geq r_i\}) = \psi(r_i) = \alpha_i+\Delta_i$.

 Ad \ref{sstep:psi-nonrightcontinuous}: Let now $r^{(n)}\searrow r_i$ monotonously decreasing and define $F_n = \{z\in \Omega: \L(z) \geq r^{(n)}\}$. Now $F_n \subset F_{n+1} \subset \cdots$ and $\bigcup_n F_n = \{z\in \Omega: \L(z) > \inf_n r^{(n)}\} =\{z\in \Omega: \L(z) > r_i\}$. The measure $\mu$ is also continuous from below, i.e.\ $\lim_{n\to\infty}\psi(r^{(n)}) = \lim \mu(F_n) = \mu(\bigcup_n F_n) = \mu(\{z\in \Omega: \L(z) > r_i\} = \alpha_i$ which is strictly less than $\alpha_i + \Delta_i = \psi(r_i)$.
 
 Ad \ref{sstep:psi-i-surjective}: Regardless of a specific case 0, A, B, or AB and regardless of which $i$ we consider, the mapping $\psi_i$ needs to be continuous on the interior of $I_i$, or else the discrete part of the measure $\Lmu$ would have an additional contribution (but we assumed the collection $(r_i)_i$ of atoms of $\Lmu$ to be complete\footnote{In the sense of ``not missing an entry'', not in the sense of completeness with respect to convergence.}).
 
 We begin with the case that $i=1,\ldots,N-1$ (where we do not have to take care of whether we are in case 0, A, B, or AB). If $\alpha_{i+1}+\Delta_{i+1} = \alpha_i$, then there actually is nothing to show as $\psi_i$ is constant. Consider then $\alpha_{i+1}+\Delta_{i+1} < \alpha_i$. Continuity at the right boundary of $I_i$, i.e.\ at $r=r_{i+1}$ follows from left-continuity of $\psi$ at this point (step~\ref{step:mappingPsi}.\ref{sstep:psi-leftcontinuous}). Monotonicity is inherited from step~\ref{step:mappingPsi}.\ref{sstep:psi-nonincreasing} and it remains to prove surjectivity. But this is clear from the construction of the intervals $I_i$ and $J_i$: It holds that $\lim_{r\searrow r_i}\psi_i(r) = \alpha_i$ (step~\ref{step:mappingPsi}.\ref{sstep:psi-nonrightcontinuous}) and $\psi(r_{i+1}) = \alpha_{i+1}+\Delta_{i+1}$. Hence by continuity, the map $\psi_i$ is surjective. 
 
 In the case that $i=0$, the first subcase (i.e.\ $\psi_0: \{\inf\L\}\to \{1\})$ is again trivial and there is nothing we need to prove. We consider the case where $\inf \L < r_1$ and $\alpha_1+\Delta_1 < 1$. We see that $\psi(\inf\L) = \mu(\{z\in \Omega: \L(z)\geq \inf \L\}) = 1$, also $\psi(r_1) = \alpha_1 + \Delta_1$. Again by continuity, $\psi_0$ needs to be surjective.
 
 The case $i=N$ is shown ``anti-parallely'' to $i=0$.
 Ad \ref{sstep:phi}: Follows directly from definition.
\end{proof}

\begin{step} \label{step:definition-ralpha}
Set $\alpha\in[0,1]$ and define $r_\alpha = \inf\{r\in \R: ~ \psi(r) \leq \alpha\}$. Then:
    \begin{enumerate}
    \item\label{sstep:ralpha_r<ralpha_psi(r)>alpha}
        For any $r < r_\alpha$ we have $\psi(r) > \alpha$. 
    \item\label{sstep:ralpha-alphaInAi}
        If $\alpha \in A_i = [\alpha_i, \alpha_i + \Delta_i)$ for some $i$, then $r_\alpha = r_i$ and $\psi(r_\alpha) = \psi(r_i) = \alpha_i + \Delta_i$.
    \item\label{sstep:ralpha-alphaNotInAi}
        If $\alpha\not\in  \bigcup_i [\alpha_i, \alpha_i + \Delta_i)$, then $\psi(r_\alpha) = \alpha$.
    \item\label{sstep:ralpha-nonincreasing}
        $\alpha\mapsto r_\alpha$ is non-increasing: For any $\alpha\geq \beta \geq 0$, we have $r_\alpha \leq r_\beta$.
    \item\label{sstep:ralpha-ralphaGreaterri-alphaSmalleralphai}
        If $r_\alpha > r_i$, then $\alpha < \alpha_i$, or equivalently, $\alpha\geq \alpha_i$, then $r_\alpha \leq r_i$.
    \item\label{sstep:alternativecharacterization}
        $r_\alpha = \sup\{r\in \R: \psi(r) > \alpha\}$.
    \item\label{sstep:tildeL}
        If we define $\tilde r_\alpha = \inf\{r\in \R: \phi(r) \leq \alpha\} = \sup\{r\in\R: \phi(r) > \alpha\}$ with $\phi$ from \ref{step:mappingPsi}.\ref{sstep:phi}, then $\tilde r_\alpha = r_\alpha$. As $\phi(r) = X(r)$, we also have $\tilde r_\alpha = \tilde \L(\alpha)$ (as defined in Section~\ref{sec:mainParadigm}).
    \end{enumerate}
\end{step}

\begin{proof}[Proof of step~\ref{step:definition-ralpha}]
Ad \ref{sstep:ralpha_r<ralpha_psi(r)>alpha}: If $r < r_\alpha$, then by definition of $r_\alpha$ we have $\psi(r) > \alpha$.

Ad \ref{sstep:ralpha-alphaInAi}: First, we know that $\psi(r_i) = \alpha_i + \Delta_i > \alpha$ and for $r < r_i$, by monotonicity, $\psi(r) > \alpha$ as well. Secondly, $\lim_{r\searrow r_i} \psi(r) = \alpha_i \leq \alpha$ and thus, again by monotonicity, $\psi(r) \leq \alpha$ for any $r \geq r_i$. This proves that 
\[\{r\in \R: \psi(r) \leq \alpha\} = (r_i, \infty). \]
Taking the infimum, we obtain $r_\alpha = r_i$ and thus $\psi(r_\alpha) = \psi(r_i) = \alpha_i + \Delta_i$.

Ad \ref{sstep:ralpha-alphaNotInAi}: In this case, we are in one of the following scenarios:
\begin{itemize}
    \item $\alpha \in [\alpha_{i+1}+\Delta_{i+1}, \alpha_i) = J_i$ for $i=1,\ldots,N-1$ with $\alpha_{i+1}+\Delta_{i+1} < \alpha_i$. (The case $\alpha_{i+1}+\Delta_{i+1} = \alpha_i$ is impossible because then $\alpha = \alpha_i$ and then we are actually in the territory of \ref{step:definition-ralpha}. \ref{sstep:ralpha-alphaInAi}).
    
    By surjectivity (step~\ref{step:mappingPsi}.\ref{sstep:psi-i-surjective}) of $\psi|_{I_i}$ onto $J_i$, the set $\{r\in I_i: \psi(r) = \alpha\}$ is non-empty. Now for any $r'\in \{r\in I_i: \psi(r) < \alpha\}$ and $r'' \in \{r\in I_i: \psi(r) = \alpha\}$, we have $r' \geq r''$. This means that $r_\alpha =\inf\{r\in I_i: \psi(r)\leq \alpha\} = \inf \{r\in I_i: \psi(r) = \alpha\}$. In particular, $r_\alpha\in \overline{I_i}$ with $\overline{I_i}$ being the closure of $I_i$ and there is a minimizing sequence $r^{(n)}\searrow r_\alpha \in I_i$  with $\psi(r^{(n)}) = \alpha$.

Now we claim that $r_\alpha > r_i$. This is true because assuming $r_\alpha = r_i$, then $\psi(r^{(n)}) = \alpha \xrightarrow{n\to\infty} \alpha$ (as the trivial limit of a constant sequence) but by step~\ref{step:mappingPsi}.\ref{sstep:psi-nonrightcontinuous} we know that $\psi(r^{(n)})\to \alpha_i > \alpha$, which is impossible. Hence, $r_\alpha \in I_i$ and $\psi(r_\alpha) = \alpha$.

\end{itemize}

\begin{itemize}
    \item Case 0, A, i.e.\ $r_N < \sup \L$ and $\alpha\in [0, \alpha_N) = J_N$.
    
    This is proven in the same way as the previous case.
    \item Case 0, B, i.e.\ $r_1 > \inf \L$ and $\alpha\in [\alpha_1 + \Delta_1, 1] = J_0$
    
    Even easier because there is no ``forbidden'' right boundary in the domain $J_0$.
    \item Case A, AB, i.e.\ $r_1 = \inf \L$ and $\alpha \in \{1\} = J_0$, which means that $\alpha = 1$.
    
    Again, this is trivial because in this case obviously $\{r\in I_0 = \{\inf \L\}: \psi(r)\leq \alpha = 1\} = \{\inf\L\}$ and the infimum of this set is its one and only element, hence $r_\alpha = r_1 = \inf \L$ and by direct evaluation, $\psi(r_\alpha) = \alpha_1+\Delta_1 = 1$. 
\end{itemize}

Ad \ref{sstep:ralpha-nonincreasing}: Follows immediately from definition.

Ad \ref{sstep:ralpha-ralphaGreaterri-alphaSmalleralphai}: We prove the equivalent form: Let $\alpha \geq \alpha_i$. We know that $r_{\alpha_i} = r_i$ (from step~\ref{step:definition-ralpha}.\ref{sstep:ralpha-alphaInAi}) and thus by setting $\beta=\alpha_i$ in step~\ref{step:definition-ralpha}.\ref{sstep:ralpha-nonincreasing}, we get $r_\alpha \leq r_{\alpha_i} = r_i$. 

Note that this is not just a corollary from step~\ref{step:definition-ralpha}.\ref{sstep:ralpha_r<ralpha_psi(r)>alpha}, but a stronger statement in the case that $r=r_i$: Step~\ref{step:definition-ralpha}.\ref{sstep:ralpha_r<ralpha_psi(r)>alpha} combined with \ref{step:definition-ralpha}.\ref{sstep:ralpha-alphaInAi} just yields the weaker statement $\alpha < \psi(r_i) = \alpha_i + \Delta_i$.

Ad \ref{sstep:alternativecharacterization}: This is elementary.

Ad \ref{sstep:tildeL}: This is due to the fact that $\phi$ is the right-continuous version of $\psi$. Take the supremization sequence $r_n$ in the alternative definition of $r_\alpha$ in \ref{step:definition-ralpha}.\ref{sstep:alternativecharacterization}, i.e.\ $\psi(r_n) > \alpha$, $r_n \nearrow r_\alpha$ and $\psi(r) \leq \alpha$ for any $r > r_\alpha$. We can choose a subsequence such that $r_n \neq r_i$ for any $i$. Then $\phi(r_n) = \psi(r_n) > \alpha$ and $r_n \nearrow r_\alpha$. Because $\phi(r)$ is less or equal than $\psi$, we also have $\psi(r) \leq \alpha$ for any $r > r_\alpha$. This means that $r_n$ is also a supremization sequence for $\tilde r_\alpha$.
\end{proof}

\begin{step} \label{step:equivalenceOfSets}
Set $\alpha\in[0,1]$. Then\footnote{The reader might wonder why we defined $\psi$ as we did, i.e.\ with an $\geq$ instead of an $>$: The central identity $\mu(\{z\in \Omega: \L(z) > \L(x)\}) \leq \alpha$ carries an $>$, after all! The proof indeed also works with a version of $\psi$ using the $>$ sign but the reasoning gets a lot more complicated at a later point: Rather than applying a left-continuous function $\psi$ to a non-increasing right-continuous function $\alpha\mapsto r_\alpha$, which yields a right-continuous function; We would have two right-continuous functions and non-increasing functions. In concatenation with each other, this would yield a function with no clear continuity property, which is unfortunate. In short, by defining $\psi$ like this, we obtain a concatenation $\alpha\mapsto \psi(r_\alpha)$ with ``nice'' properties with the drawback that step~\ref{step:equivalenceOfSets} is slightly more involved.}
$\{x\in \Omega: \mu(\{z\in \Omega: \L(z) > \L(x)\}) \leq \alpha\} = \{x\in \Omega: \L(x) \geq r_{\alpha}\}.$
\end{step}

\begin{proof}[Proof of step~\ref{step:equivalenceOfSets}]
``$\subseteq$'': Take $x\in \Omega$ such that $\mu(\{z\in \Omega: \L(z) > \L(x)\}) \leq \alpha$. Assume that $\L(x) < r_\alpha$. Then by Step~\ref{step:definition-ralpha}.\ref{sstep:ralpha_r<ralpha_psi(r)>alpha}, $\alpha< \psi(\L(x)) = \mu(\{z\in \Omega: \L(z) \geq \L(x)\}) = \mu(\{z\in \Omega: \L(z) > \L(x)\}) + \mu(\{z\in \Omega: \L(z) = \L(x)\})\leq \alpha + \mu(\{z\in \Omega: \L(z) = \L(x)\})$, hence 
\[\mu(\{z\in \Omega: \L(z) = \L(x)\}) > 0.\]
This is only possible if $\L(x) = r_i$ with $r_i = \L(x) < r_\alpha$ by assumption. The last statement yields (using Step~\ref{step:mappingPsi}.\ref{sstep:ralpha-ralphaGreaterri-alphaSmalleralphai}) that $\alpha < \alpha_i$. But this means that 
\[\alpha \stackrel{\text{assump.}}{\geq}\mu(\{z\in \Omega: \L(z) > \L(x)\}) = \mu(\{z\in \Omega: \L(z) > r_i\}) = \alpha_i \stackrel{\text{step~\ref{step:mappingPsi}.\ref{sstep:ralpha-ralphaGreaterri-alphaSmalleralphai}}}> \alpha,  \]
a contradiction.

``$\supseteq$'': Let $x\in \Omega$ such that $\L(x)\geq r_\alpha$. We need to distinguish two cases. 

Case I: $\L(x) = r_\alpha$.  Then 
\begin{align*}
    \mu(\{z\in \Omega: \L(z) > \L(x)\}) &= \mu(\{z\in \Omega: \L(z) > r_\alpha\})\\
    &= \mu(\{z\in \Omega: \L(z) \geq  r_\alpha\}) - \mu(\{z\in \Omega: \L(z) = r_\alpha\})\\
    &= \psi(r_\alpha) - \mu(\{z\in \Omega: \L(z) = r_\alpha\})\\
    &=\begin{cases}
    (\alpha_i + \Delta_i) - \Delta_i \leq \alpha, &\text{ for } \alpha\in [\alpha_i, \alpha_i+\Delta_i) \text{ for some } i\\
    \alpha - ``\geq 0'' \leq \alpha, &\text{ else}
    \end{cases}
\end{align*}

Case II: $\L(x) > r_\alpha$. 

Then by definition of $r_\alpha$, $\alpha \geq \psi(\L(x)) = \mu(\{z\in \Omega: \L(z) \geq \L(x)\}) \geq \mu(\{z\in \Omega: \L(z) > \L(x)\})$.

This means that in each case, $\mu(\{z\in \Omega: \L(z) > \L(x)\})\})\leq \alpha$ which concludes the proof of inclusion.
 
\end{proof}

\begin{step}\label{step:connection-muAndAlphaalphai}
If $\alpha \not \in \bigcup_i [\alpha_i, \alpha_i + \Delta_i)$, then $\mu\left(\left\{x\in \Omega:\mu\left(\left\{z\in \Omega: \L(z) > \L(x)\right\}\right)\leq \alpha\right\}\right) = \alpha$.

If $\alpha \in \bigcup_i [\alpha_i, \alpha_i + \Delta_i)$, then $\mu\left(\left\{x\in \Omega:\mu\left(\left\{z\in \Omega: \L(z) > \L(x)\right\}\right)\leq \alpha\right\}\right) = \alpha_i + \Delta_i$.
\end{step}

\begin{proof}[Proof of step~\ref{step:connection-muAndAlphaalphai}] We calculate (for the first equality we use step~\ref{step:equivalenceOfSets})
\begin{align*}
    &\mu\left(\left\{x\in \Omega:\mu\left(\left\{z\in \Omega: \L(z) > \L(x)\right\}\right)\leq \alpha\right\}\right)\\
    \qquad&= \mu(\{x \in \Omega: \L(x)\geq r_\alpha\})\\
    \qquad&= \psi(r_\alpha)\\
    \qquad&= \begin{cases}\alpha,&\text{ if }\alpha \not \in \bigcup_i [\alpha_i, \alpha_i + \Delta_i) \\ \alpha_i + \Delta_i,&\text{ if }\alpha \in \bigcup_i [\alpha_i, \alpha_i + \Delta_i)\end{cases}
\end{align*}
 
\end{proof}
This concludes the proof of Lemma~\ref{lem:main}.
\end{proof}

\end{lemma}

\begin{lemma}\label{prop:push-forward}
	Let $\Omega \subseteq \mathbb{R}^d$ with $d\in\mathbb{N}$
	and 
	$\L: (\Omega, \mathcal{B}^d) \to (\mathbb{R}, \mathcal{B})$ with Borel $\sigma$-algebra $\mathcal{B}^d$ on $\Omega$ and $\mathcal{B}$ on $\mathbb{R}$. 
	Let further $(\Omega, \mathcal{B}^d, \mu)$ be a probability space. We assume that we know the discrete part of $\Lmu$ as in Lemma~\ref{lem:main}, i.e.\
	\begin{align*}
    \Lmu(\{r_i\}) &= \mu(\{z\in \Omega: \L(z) = r_i\}) = \Delta_i\\
    \Lmu((r_i, \infty)) &=  \mu(\{z\in \Omega: \L(z) > r_i\}) = \alpha_i.
\end{align*}
    We assume the $\alpha_i$ to be ordered, i.e.\ $\alpha_1 < \alpha_1 + \Delta_1 \leq \alpha_2 < \cdots$. If there are only $N$-many jumps, we set $\alpha_{N+1} = 1$.
	
	Then, $\Phi: \Omega \to [0,1]$ given by 
	\begin{equation}
		\Phi(x) = X(\L(x)) = \int_{\L(z) > \L(x)} \mathrm d\mu(z) \label{eq:Xx}
	\end{equation}
	is distributed according to the law
	\[ \sum_{i} \Delta_i\cdot \delta_{\alpha_i} + \operatorname{Unif} \left(\left[0,1\right] \setminus \bigcup_i [\alpha_i, \alpha_i+\Delta_i]\right),\]
	where $\delta_{p}$ denotes the Dirac measure at $p$ and  $\operatorname{Unif}(I)$ is a uniform distribution on the set $I$.
\end{lemma}
\begin{proof}
The law of $\Phi\in [0,1]$ is the cumulative distribution function of the push-forward $\Phi_\#\mu$ given by 

\begin{align*}
    \mu\circ \Phi^{-1}([0, \alpha]) &=  \mu \left(\left\lbrace x \in \Omega : \Phi(x) \in [0,\alpha] \right \rbrace \right)\\
    &=\mu \left(\left \lbrace x \in \Omega : \int_{\left \lbrace z\in \Omega: \L(z) > \L(x) \right \rbrace} \mathrm d \mu(z) \le \alpha \right\rbrace  \right) \\
    &= \mu \left( \left\lbrace x\in \Omega: \mu \left(\left\lbrace z\in \Omega: \L(z) > \L(x) \right\rbrace \right) \le \alpha \right\rbrace \right)\\
    &= \begin{cases} \alpha_i + \Delta_i &\text{for } \alpha \in [\alpha_i, \alpha_i + \Delta_i)\\ \alpha &\text{else.}
\end{cases}
\end{align*}
This is exactly the cumulative distribution function of the measure 
$ \sum_{i} \Delta_i\cdot \delta_{\alpha_i} + \operatorname{Unif}\left(\left[0,1\right] \setminus \bigcup_i [\alpha_i, \alpha_i+\Delta_i]\right).$
\end{proof}

\begin{corollary}\label{cor:inversion}
	Let $\Omega \subseteq \mathbb{R}^d$ with $d\in\mathbb{N}$
	and 
	$\L: (\Omega, \mathcal{B}^d) \to (\mathbb{R}, \mathcal{B})$ with Borel $\sigma$-algebra $\mathcal{B}^d$ on $\Omega$ and $\mathcal{B}$ on $\mathbb{R}$. We assume $\L$ is bounded from above and below and that $\mu\circ\L^{-1}$ has no discrete atoms, i.e.\ for any $r\in [\inf\L, \sup\L]$, we have $\mu(\{x\in \Omega: \L(x) = r\}) = 0$.
	
	Then if $x\in\Omega$ is a random variable distributed according to $\mu$, then $\Phi(x)$ is uniformly distributed on $[0,1]$.
\end{corollary}
\begin{proof}
Follows immediately from Proposition~\ref{prop:push-forward}.
\end{proof}

\begin{lemma}\label{lem:integral_wrt_pshfwd}
    With the assumptions and definitions of Proposition~\ref{prop:push-forward}, we have that 
    \begin{equation}
        \int_0^1 \tilde \L(t) \diff \Phi_\#\mu(t) = \int_0^1 \tilde \L(t) \diff t.
    \end{equation}
\end{lemma}
\begin{proof}
First we recall from Proposition~\ref{prop:push-forward} that 
$\Phi_\#\mu = \sum_{i} \Delta_i\cdot \delta_{\alpha_i} + \operatorname{Unif}\left(\left[0,1\right] \setminus \bigcup_i [\alpha_i, \alpha_i+\Delta_i]\right).$
Thus, 
\begin{figure}
    \centering
    \includegraphics[width=1\textwidth]{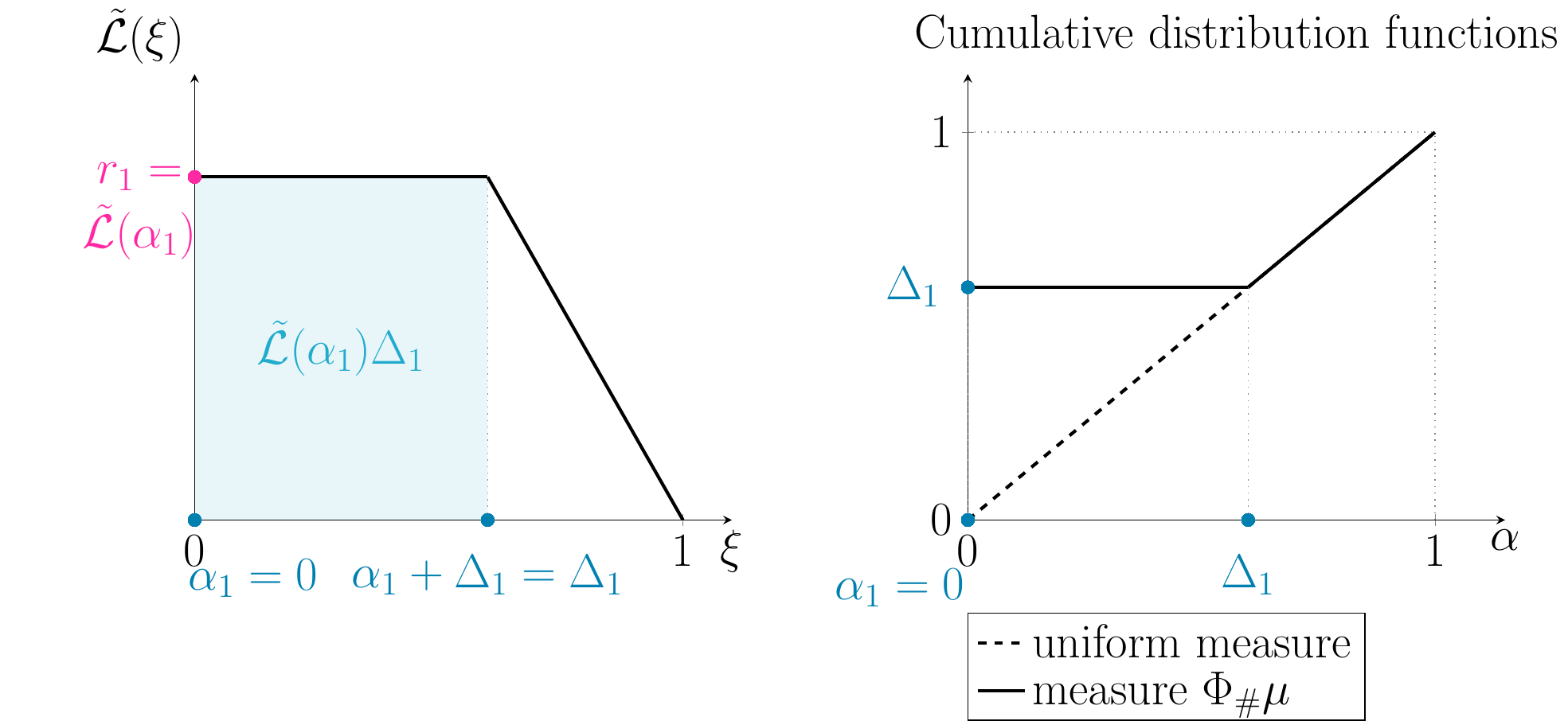}
    \caption{$\tilde \L(\alpha_i)\cdot \Delta_i = \int_{\alpha_i}^{\alpha_i+\Delta_i} \tilde \L(t) \diff t$.}
    \label{fig:UniformVSPhiHashMu}
\end{figure}
\begin{align*}
    \int_0^1 \tilde \L(t) [\diff\Phi_\#\mu(t) - \diff t] &= \sum_i \left[ \tilde \L(\alpha_i)\cdot \Delta_i -\int_{\alpha_i}^{\alpha_i+\Delta_i}\tilde \L(t)\diff t\right]\\
    &= \sum_i \int_{\alpha_i}^{\alpha_i+\Delta_i}\left[\tilde \L(\alpha_i) - \tilde \L(t)\right] \diff t
\end{align*}
By characterization \ref{step:definition-ralpha}.\ref{sstep:tildeL} in the proof of Lemma~\ref{lem:main}, we know that for $t\in [\alpha_i, \alpha_i+\Delta_i)$, 
\[\tilde \L(t) = \tilde r_t = r_t = r_i = \tilde \L(\alpha_i).\]
This means that the integrand is always zero and thus
\[\int_0^1 \tilde \L(r) [\diff\Phi_\#\mu(r) - \diff r] = 0.\]
For a visualization of this fact see Fig.~\ref{fig:UniformVSPhiHashMu}
\end{proof}

We have shown that although $\Phi_\#\mu$ is not the uniform measure on $[0,1]$, integration of $\tilde \L$ with respect to both measures yields the same value. Now we need to validate the second dubious equality, i.e.\ \eqref{eq:line1}. 
\section{Relating the likelihood with the generalized inverse of $X$}\label{sec:inverse}

We defined the generalized inverse for $X$ in Section~\ref{sec:mainParadigm} as $\tilde \L(\xi) = \sup\{\lambda\in  \im \L: X(\lambda) > \xi\}$. There, we already indicated that neither $\tilde{\L}(X(\lambda)) = \lambda$ for $\lambda \in \mathbb{R}$ holds nor does $\tilde{\L}(X(\L(x))) = \L(x)$ for all $x \in \Omega$. This is demonstrated in the following with specific examples. We then continue to show that at least $\tilde{\L}(X(\L(x))) = \L(x)$ for $\mu$-almost-all $x$. Although more restrictive, this meaning of ``$\tilde \L$ is the inverse of $X$'' is all we need in order for the equality in \eqref{eq:line1}
\begin{equation*}
    \int_\Omega \L(x) \diff \mu(x) = \int_\Omega \tilde{\L}(X(\L(x))) \diff \mu(x)
\end{equation*}
to hold. 
\subsection{$\tilde \L$ is not the inverse of $X$}
\begin{figure}
    \centering
    \includegraphics[width=1\textwidth]{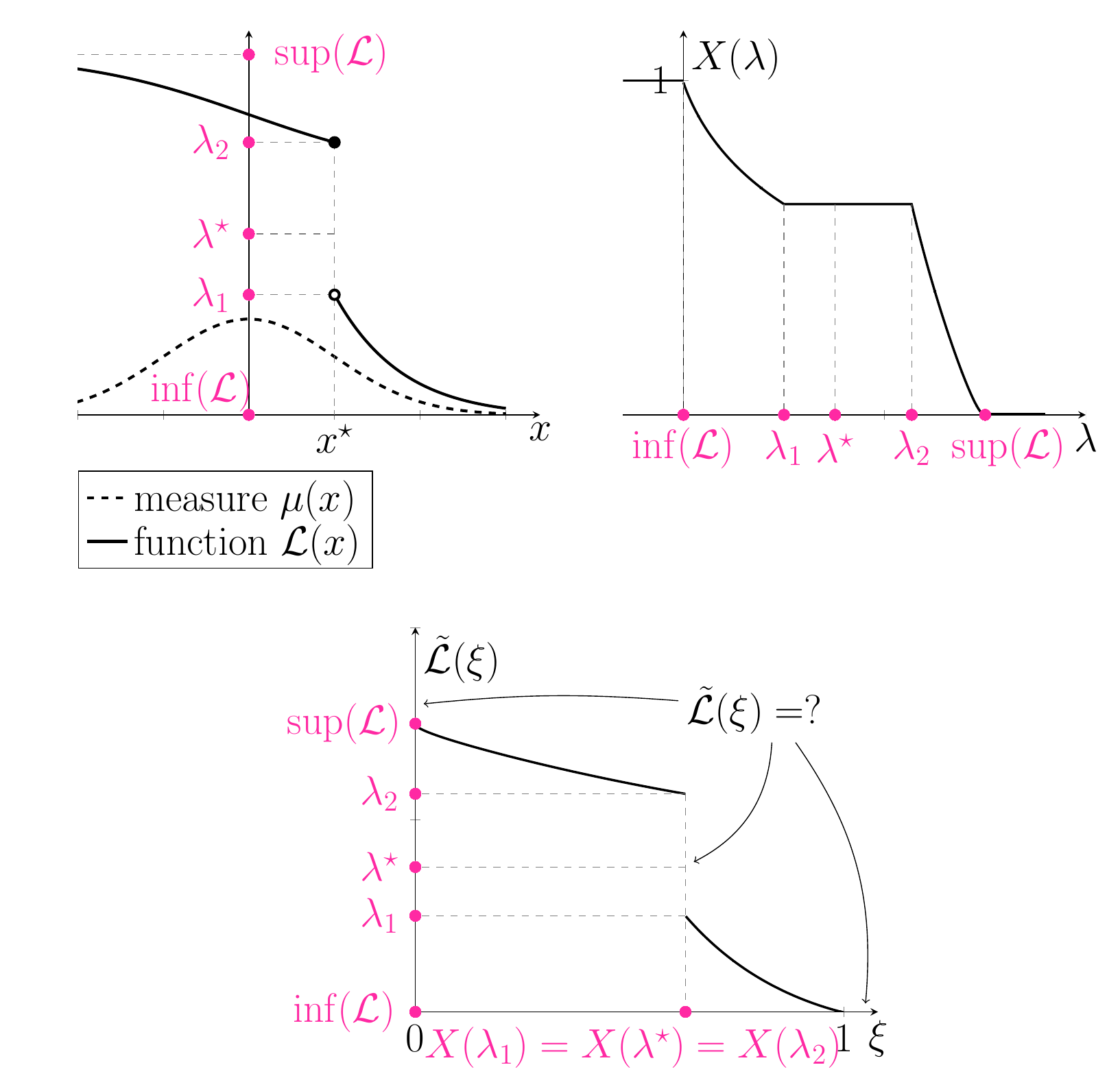} 
    \caption{The left plot shows one suitable problem setting such that $X(\lambda)$ has plateaus (right). Therefore, its inverse is not even defined especially at the jumps where $X(\lambda)$ is non-injective in these ranges (row 2).}
    \label{fig:LofXlambda-notWorking}
\end{figure}

It can be easily demonstrated that the identity $\tilde \L(X(\lambda)) \neq \lambda$ does not hold for $\lambda \in \mathbb{R}$ in general. In particular, if $X$ has a plateau, this fails: Given $\lambda_1\neq \lambda_2$ with $X(\lambda_1) = X(\lambda_2)$, there is no way to define an inverse $\tilde \L$ because $X$ is not injective. See Fig.~\ref{fig:LofXlambda-notWorking} for a concrete example of a likelihood function $\L$ such that $X$ is not injective. As is apparent, jumps in the range of $\L$ are a problem (if they are in the support of $\mu$).

The counterexample seems to suggest that only $\lambda \not \in \operatorname{range}\L$ are a problem and that we could at least expect the identity $\tilde \L \circ X = \operatorname{id}$ to hold on the range of $\L$, i.e.\ $\tilde \L(X(\L(x))) = \L(x)$, but as the next section reveals, this is not true, either.

We think about when and how $\tilde \L(X(\L(x))) = \L(x)$ is true:
\begin{align*}
    \tilde \L(X(\L(x))) &= \sup\{\lambda\in \im \L: ~X(\lambda) > X(\L(x))\}\\
    \intertext{As for $\lambda \geq \L(x)$ we have $X(\lambda) \leq X(\L(x))$, we can ignore any $\lambda \geq \L(x)$ in the specification of the set on the right hand side, hence}
    &=\sup\{\lambda\in \im \L: ~\lambda < \L(x) \quad \text{ and } 
    \quad X(\lambda) > X(\L(x))\}\\
    \intertext{By definition, $X(\lambda) = \mu(\{z\in\Omega:~ \L(z) > \lambda\})$ and thus $X(\lambda) - X(\L(x)) = \mu(\{z\in\Omega:~ \L(x) \geq \L(z) > \lambda\})$ (this is due to the elementary fact that for sets $A\subset B$ and a measure $\nu$, we have $\nu(B) - \nu(A) = \nu(B\setminus A)$).}
    &=\sup\Big\{\lambda\in \im \L: ~\lambda < \L(x) \quad \text{ and } 
    \quad\mu\big(\{z\in\Omega: \L(z)\in (\lambda, \L(x)]\}\big) > 0\Big\}
\end{align*}
We need to discern two cases now:

\textbf{Case I: For all $\lambda \in \im \L$ with $\lambda < \L(x)$, we have that $\mu\big(\{z\in\Omega: \L(z)\in (\lambda, \L(x)]\}\big) > 0$.}
As for $\lambda = \L(x)$ and thus $(\lambda, \L(x)] = \emptyset$, we clearly have $\mu\big(\{z\in\Omega: \L(z)\in (\lambda, \L(x)]\}\big) = 0$, we know that 
\begin{align*}
&\sup\Big\{\lambda\in \im \L: ~\lambda < \L(x) \text{ and } \mu\big(\{z\in\Omega: \L(z)\in (\lambda, \L(x)]\}\big) > 0\Big\}\\
&\quad= \sup\{\lambda \in \im \L: \lambda < \L(x)\}\\ &\quad= \L(x),    
\end{align*}
i.e.\ indeed 
\[\tilde \L(X(\L(x))) = \L(x).\]

\begin{figure}
    \centering
    \includegraphics[width=1\textwidth]{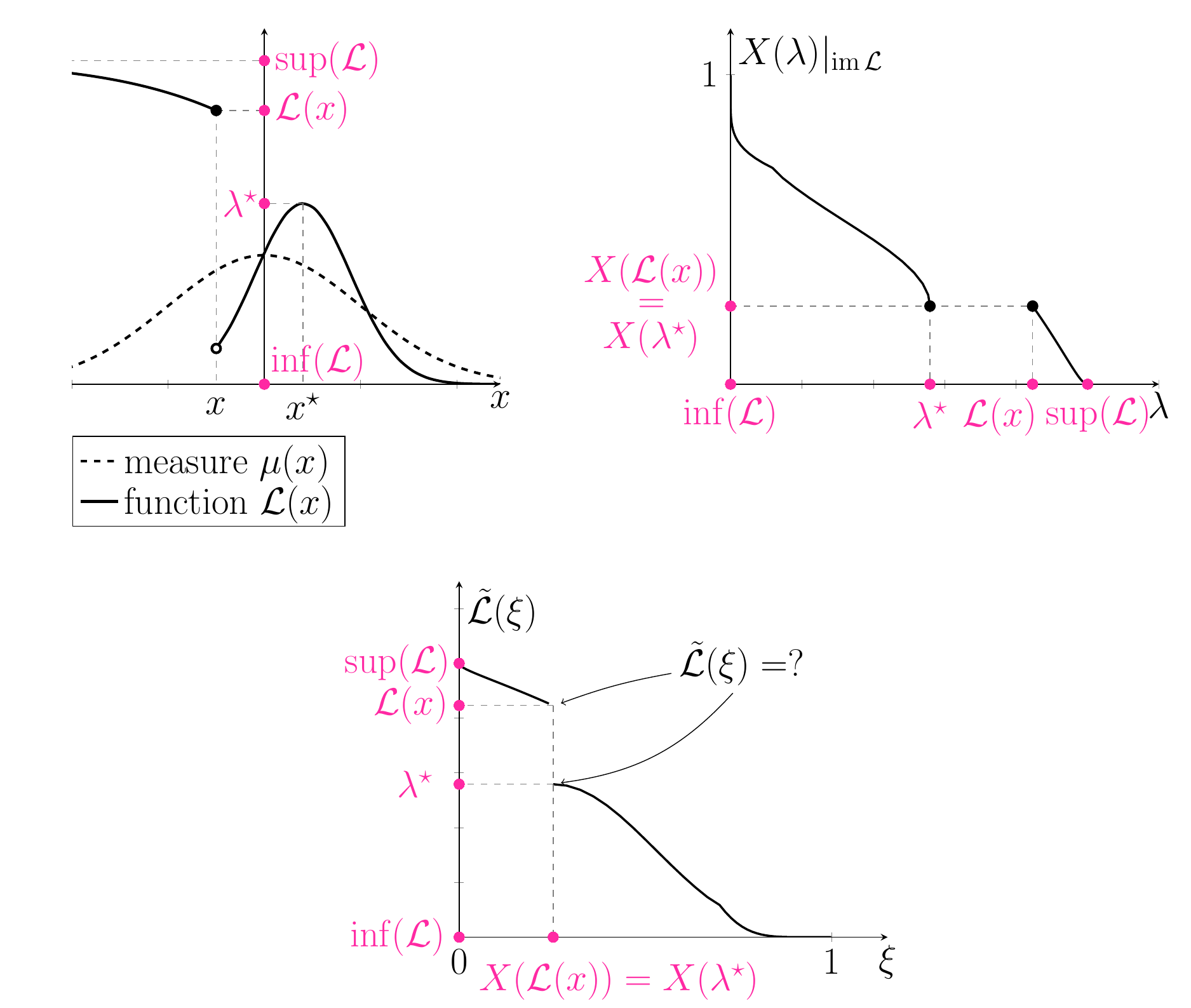}
    \caption{In the problem setting shown in the left plot, the property $\tilde\L(X(\L(x))) = \L(x)$ does \textit{not} hold because $x$ and $x^\star$ result in the same value for $X$ (right) and, therefore, $\tilde\L(\xi)$ is not defined especially at this point (row 2).}
    \label{fig:counterexample-case2}
\end{figure}

\textbf{Case II: (Complement of Case I)}
There exists a $\lambda^\star\in \im \L$, i.e.\ $\lambda^\star = \L(x^\star), x^\star \in \Omega,$ such that $\lambda < \L(x)$ with the property that
\[ \mu\big(\{z\in\Omega: \L(z)\in (\lambda^\star, \L(x)]\}\big) = 0,\]

then of course 
\begin{align*}
&\sup\Big\{\lambda\in \im \L: ~\lambda < \L(x) \text{ and } \mu\big(\{z\in\Omega: \L(z)\in (\lambda, \L(x)]\}\big) > 0\Big\} \leq \lambda^\star
\end{align*}
and thus 
\[\tilde \L(X(\L(x))) \leq \lambda^\star = \L(x^\star) < \L(x),\]
i.e.\ the desirable property $\tilde\L(X(\L(x))) = \L(x)$ does \textit{not} hold here.
It is not difficult to construct a concrete example for this case, see Fig.~\ref{fig:counterexample-case2}.

\subsection{$\tilde \L$ is almost everywhere the inverse of $X$ on its range}\label{sec:inverseX}
We saw that the $\tilde \L(X(\L(x))) = \L(x)$ does not hold for all $x\in \Omega$ but that it is violated if for a given $x\in \Omega$, there exists a $\lambda^\star = \L(x^\star) < \L(x)$ such that 
\[ \mu\big(\{z\in\Omega: \L(z)\in (\lambda^\star, \L(x)]\}\big) = 0.\]
We will call such points $x \in \Omega$ \textit{problematic points} in this context. We are interested in whether the set of \textit{problematic points} has vanishing $\mu$-measure. If that is the case, then the identity $\tilde \L(X(\L(x))) = \L(x)$ still holds $\mu$-almost-everywhere.
The following Lemma shows that this is the case
\begin{lemma}
    \label{lem:characterization-problematic}
    The set of problematic points has $\mu$-measure $0$, i.e. for $\mu$-almost-all $x\in\Omega$, 
    \[\tilde \L(X(\L(x))) = \L(x). \]
    \end{lemma}
\begin{proof}
The set of problematic points is defined as 
\begin{align*}P&=\Big\{x\in\Omega: \text{ There is a } \lambda^\star = \L(x^\star) < \L(x) \\  &\qquad\text{ such that } \mu(\{z\in\Omega:~ \L(z)\in(\lambda^\star, \L(x)]\}) = 0 \Big\}
\end{align*}
Every problematic point $x_p\in P$ ``generates'' an interval\footnote{Here, $p$ is an index which is not necessarily countable but which indexes the whole family of points $x_p \in P$.} $I_p = (\lambda_p^\star, \L(x_p)]$ which has vanishing $\L_\#\mu$-measure, i.e.\ $\mu(\{z\in \Omega: \L(z) \in (\lambda_p^\star, \L(x_p)]\}) = 0$, or $\L_\#\mu((\lambda_p^\star, \L(x_p)]) = 0$. The (possibly uncountable) union of those intervals $I_p$ over all $x_p$ becomes a \textit{countable} union of intervals,
\[\bigcup_p I_p = \bigcup_{n\in \N} (a_n, b_n] \cup \bigcup_{m\in\N} (c_m, d_m),\]
where $a_n$ and $c_m$ are limits of certain subsequences $\lambda_{p_m}^{(n)}$ and $b_n$ and $d_m$ are limits of certain subsequences $\L(x_{p_m}^{(n)})$, respectively. This follows from Lemma~\ref{lem:uncountable-union-halfopen}. By continuity from below for $\L_\#\mu$, we also have that $\L_\#\mu((a_n, b_n]) = \mu(\{z\in\Omega: \L(z) \in (a_n, b_n]\}) = 0$, for any $n\in \N$, and the same for the intervals $(c_m, d_m)$. As countable sums of terms with value $0$ are $0$, we obtain $\L_\#\mu(\bigcup_p I_p) = 0$.

Now note that for every $x_p\in P$, by construction, $\L(x_p) \in I_p$, which in turn means that
\[ P \subseteq \left\{x\in\Omega: ~ \L(x) \in \bigcup_p I_p\right\} = \L^{-1}\left(\bigcup_p I_p\right). \]
Finally, we obtain $\mu(P) \leq \mu \circ \L^{-1}\left(\bigcup_p I_p\right) = \L_\#\mu\left(\bigcup_p I_p\right) = 0$.
\end{proof}


\section{Measure-theoretical statements}\label{appendix-MeasureTheoreticalStatements}
    We recall the following fact from measure theory.
    \begin{lemma}\label{lem:countableunion}
    Any open subset $U\subseteq \R$ can be written as an at most countable union of disjoint open intervals, i.e.\ \[ U = \bigcup_{n\in\N} (c_n, d_n)\]
    with $d_n > c_n$ and $(c_n, d_n) \cap (c_m, d_m) = \emptyset$ for $m\neq n$.
    \end{lemma}
    
    \begin{lemma} \label{lem:uncountable-union-open}
    Let $\Pi$ be a set (possibly uncountable) of indices and $\{I_p\}_{p\in \Pi}$ be a (possibly uncountable) family of intervals of the form $I_p = (a_p, b_p)$ with $b_p > a_p$ for any index $p\in \Pi$. Then 
    \begin{equation}\bigcup_{p\in\Pi} I_p = \bigcup_{n\in\N} (c_n, d_n) \label{eq:openunion}
    \end{equation} 
    with $d_n > c_n$ and  $(c_n, d_n) \cap (c_m, d_m) = \emptyset$ for $m\neq n$, i.e.\ the uncountable union simplifies to a countable union of open intervals. 
    
    Furthermore, for any $n\in\N$, there is a sequence of indices $\{p_m^{(1)}\}_m \subset \Pi$ such that $c_n = \lim_{m\to \infty} a_{p_m^{(1)}}$ and a sequence of indices $\{p_m^{(2)}\}_m \subset \Pi$ such that $d_n = \lim_{m\to \infty} b_{p_m^{(2)}}$.
    \end{lemma}
    \begin{proof}
    An arbitrary union (even uncountable) of open sets is itself open, thus $U =  \bigcup_{p\in\Pi} I_p$ is open. Then we can apply Lemma~\ref{lem:countableunion} to $U$ which immediately yields the expression \eqref{eq:openunion}. Now we define $\Pi_i$ as the subset of indices such that 
    \[ \bigcup_{p\in \Pi_i} I_p = (c_i, d_i).\]
    This is possible because each interval $I_p$ can only be a subset of exactly one of the disjoint sets $(c_n, d_n)$. Now clearly $c_i = \inf_{p\in \Pi_i} a_p$ and $d_i = \sup_{p\in \Pi_i} b_p$. Then $\{p_m^{(1)}\}_m$ can be chosen as the infimization sequence for $c_i$ and analogously for $d_i$.
    \end{proof}

    \begin{lemma} \label{lem:uncountable-union-halfopen}
    Let $\{I_p\}_{p\in \Pi}$ be a (possibly uncountable) family of intervals of the form $I_p = (a_p, b_p]$ with $b_p > a_p$ for any $p\in \Pi$. Then 
    \[ \bigcup_{p\in\Pi} I_p = \bigcup_{n\in\N} J_n\]
    where each $J_n$ is either $J_n = (c_n, d_n)$ or $(c_n, d_n]$ with $c_n, d_n$ being the same numbers as in Lemma~\ref{lem:uncountable-union-open}, i.e.\ the uncountable union simplifies to a countable union of open and half-open intervals.
    \begin{proof}
    From the proof of Lemma~\ref{lem:uncountable-union-open}, we know that 
    \[\bigcup_{p\in \Pi_i} I_p = (c_n, d_n) \cup \bigcup_{p \in \Pi_i}\{b_p\}.\]
    Now there are two cases to be distinguished: If there is a $p^\star\in \Pi_i$ such that $b_{p^\star} = \sup_{p\in\Pi_i} b_p = d_n,$ then $\bigcup_{p\in \Pi_i} I_p = (c_n, b_{p^\star}] = (c_n, d_n]$. If that is not the case, i.e.\ if for any $p\in \Pi_i$, we have $b_p < d_n$, then $\bigcup_{p\in \Pi_i} I_p = (c_n, d_n)$.
    \end{proof}
    \end{lemma}

%
%
\end{appendix}

 \section*{Acknowledgements}
 Both authors want to thank Maria Neuss-Radu for fruitful discussion. 
    D.S. likes to acknowledge the support by the Bavarian Equal Opportunities Sponsorship -- Realisierung von Chancengleichheit von Frauen in Forschung und Lehre (FFL) -- Realization Equal Opportunities for Women in Research and Teaching.“ 
%
%



\bibliographystyle{imsart-number} 
\bibliography{paper-ref.bib}       


\end{document}